\documentclass[a4paper,12pt]{article}
\usepackage{amsmath}
\usepackage{amsfonts}
\usepackage{amssymb}
\usepackage{latexsym}
\usepackage{epsfig}
\usepackage{graphicx}
\usepackage{oldgerm}
\usepackage{theorem}

\def\C{\mathbb{C}}
\def\R{\mathbb{R}}
\def\N{\mathbb{N}}
\def\W{\mathbb{W}}
\def\Z{\mathbb{Z}}

\def\x{\mib{x}}
\def\z{\mib{z}}
\def\u{\mib{u}}
\def\v{\mib{v}}
\def\f{\mib{f}}
\def\t{\mib{t}}

\def\X{\mib{X}}
\def\bZ{\mib{Z}}

\def\S{\mib{S}}

\def\rP{{\rm P}}
\def\rE{{\rm E}}

\def\P{{\mathbb P}}
\def\E{{\mathbb E}}

\def\1{{\bf 1}}

\def\mbK{\mathbb{K}}
\def\bK{{\bf K}}
\def\bE{{\bf E}}
\def\bP{{\bf P}}

\def\mM{\mathfrak{M}}

\def\sfM{{\sf M}}
\def\sfK{{\sf K}}

\def\cF{{\cal F}}
\def\cM{{\cal M}}
\def\cD{{\cal D}}
\def\cG{{\cal G}}
\def\cB{{\cal B}}
\def\cP{{\cal P}}

\def\I{\mathbb{I}}
\def\J{\mathbb{J}}

\def\dto{\stackrel{\rm d}{\to}}
\def\vXi{{\bf \Xi}}

\theorembodyfont{\itshape}

\newtheorem{thm}{Theorem}[section]
\newtheorem{lem}[thm]{Lemma}

\newtheorem{prop}[thm]{Proposition}

\newtheorem{df}[thm]{Definition}

\newcommand{\mib}[1]{\mbox{\boldmath $#1$}}
\newcommand{\SSC}[1]{\section{#1}\setcounter{equation}{0}}
\newcommand{\qed}{\hbox{\rule[-2pt]{3pt}{6pt}}}

\begin{document}

\title{\bf 
Determinantal Martingales and Correlations \\
of Noncolliding Random Walks}
\author{
Makoto Katori
\footnote{
Department of Physics,
Faculty of Science and Engineering,
Chuo University, 
Kasuga, Bunkyo-ku, Tokyo 112-8551, Japan;
e-mail: katori@phys.chuo-u.ac.jp
}}
\date{20 December 2014}
\pagestyle{plain}
\maketitle
\begin{abstract}
We study the noncolliding random walk (RW), which is 
a particle system of one-dimensional, simple and symmetric RWs
starting from distinct even sites and conditioned 
never to collide with each other.
When the number of particles is finite, $N < \infty$, 
this discrete process is constructed as an $h$-transform
of absorbing RW in the $N$-dimensional Weyl chamber.
We consider Fujita's polynomial martingales of RW
with time-dependent coefficients and express
them by introducing a complex Markov process.
It is a complexification of RW, in which 
independent increments of its imaginary part are in 
the hyperbolic secant distribution,
and it gives a discrete-time conformal martingale.
The $h$-transform is represented by a determinant
of the matrix, whose entries are all polynomial martingales.
From this determinantal-martingale representation (DMR) 
of the process,
we prove that the noncolliding RW is determinantal
for any initial configuration with $N < \infty$, 
and determine the correlation kernel as a function
of initial configuration.
We show that noncolliding RWs started at infinite-particle
configurations having equidistant spacing
are well-defined as determinantal processes
and give DMRs for them.
Tracing the relaxation phenomena shown by these infinite-particle systems,
we obtain a family of equilibrium processes parameterized
by particle density, which are determinantal
with the discrete analogues of the extended sine-kernel
of Dyson's Brownian motion model with $\beta=2$.
Following Donsker's invariance principle, 
convergence of noncolliding RWs to the Dyson model 
is also discussed.
\end{abstract}



\normalsize

\SSC{Introduction \label{sec:Introduction}}

Let $\zeta$ be a random variable binomially distributed as
\begin{equation}
\rP[\zeta=1]=\frac{1}{2}, \quad
\rP[\zeta=-1]=\frac{1}{2}, 
\label{eqn:P1}
\end{equation}
so that the Laplace transform of the probability distribution is given by
\begin{equation}
\rE[e^{\alpha \zeta}]= \cosh \alpha, \quad \alpha \in \R.
\label{eqn:zeta1}
\end{equation}
For $N \in \N \equiv \{1,2, \dots \}$, let
$\{\zeta_j(t): 1 \leq j \leq N, t \in \N\}$
be a family of i.i.d.random variables 
which follow the same probability law with $\zeta$. 
We consider a random walk (RW) on $\Z^N$,  
$\S(t)=(S_1(t), \dots, S_N(t)), t \in \N_0 \equiv \{0\} \cup \N$, 
in which the components $S_j(t), j=1,2,\dots, N$ are independent
simple and symmetric RWs;
\begin{eqnarray}
S_j(0) &=& u_j \in \Z,
\nonumber\\
S_j(t) &=& u_j+\zeta_j(1)+\zeta_j(2)+ \cdots
+ \zeta_j(t), \quad
\quad t \in \N, \quad 1 \leq j \leq N.
\nonumber
\end{eqnarray}
Let $\Z_{\rm e}=2 \Z=\{\dots, -2, 0, 2, 4, \dots\}$
and $\Z_{\rm o}=1+2 \Z =\{\dots, -1,1,3,5, \dots\}$.
For each component, $S_j(\cdot), 1 \leq j \leq N$,
the transition probability is given by
\begin{eqnarray}
&& p(t-s,y|x) 
= \rP[S_j(t)=y | S_j(s)=x]
\nonumber\\
&& \quad = \left \{ \begin{array}{l}
\displaystyle{ 
\frac{1}{2^{t-s}} {t-s \choose [(t-s)+(y-x)]/2 }}, \cr
\qquad \qquad
\mbox{if $t \geq s, \, -(t-s) \leq y-x \leq t-s,  \, (t-s)+(y-x) \in \Z_{\rm e}$}, \cr
0, \quad \qquad \mbox{otherwise}.
\end{array} \right.
\label{eqn:tp1}
\end{eqnarray}
We always take the initial point 
$\u=(u_1, \dots, u_N) =\S(0)$ from $\Z^N_{\rm e}$,
then $\S(t) \in \Z^N_{\rm e}$, if $t$ is even,
and $\S(t) \in \Z^N_{\rm o}$, if $t$ is odd.
The probability space is denoted as $(\Omega, \cF, \rP_{\u})$
and expectation is written as $\rE_{\u}$.

Let
$$
\W_N=\{\x=(x_1, \dots, x_N) \in \R^N: x_1 < \cdots < x_N\}
$$
be the Weyl chamber of type A$_{N-1}$.
Define $\tau_{\u}$ to be the first exit time from the Weyl chamber
of the RW started at $\u \in \Z_{\rm e}^N \cap \W_N$, 
$$
\tau_{\u}=\inf \{ t \geq 1: \S(t) \notin \W_N\}.
$$
In the present paper, we study the RW 
{\it conditioned to stay in $\W_N$ forever}.
That is, $\tau_{\u}=\infty$ is conditioned. 
We call such a conditional RW 
the {\it (simple and symmetric) noncolliding RW}, since when we regard the
$j$-th component $S_j(\cdot)$ as the position of
$j$-th particle on $\Z, 1 \leq j \leq N$,
if $\tau_{\u} < \infty$, 
then at $t=\tau_{\u}$ there is at least
one pair of particles $(j, j+1)$, which collide with each other;
$S_j(\tau_{\u})=S_{j+1}(\tau_{\u}), 1 \leq j \leq N-1$.
Such a conditional RW is also called a system of 
{\it vicious walkers} in statistical physics \cite{Fis84,CK03},
{\it non-intersecting paths}, {\it non-intersecting walks},
and {\it ordered random walks} 
in enumerative combinatorics and probability theory (see 
\cite{Kra06,EK08} and Chapter 10 in \cite{For10}).

Let $\mM$ be the space of nonnegative integer-valued Radon measure 
on $\Z$
and $\mM_0 \equiv \{\xi \in \mM : \xi(\{x\}) \leq 1,
\forall x \in \Z\}$.
We consider the noncolliding RW as a process in $\mM_0$
and represent it by
\begin{equation}
\Xi(t, \cdot)=\sum_{j=1}^N \delta_{S^0_j(t)}(\cdot),
\quad t \in \N_0,
\label{eqn:Xi1}
\end{equation}
where 
\begin{equation}
\S^0(t)=(S^0_1(t), \dots, S^0_N(t)) \in \Z^N \cap \W_N, 
\quad t \in \N_0.
\label{eqn:Xi2}
\end{equation}
The configuration $\Xi(t, \cdot) \in \mM_0, t \in \N_0$ is
unlabeled, while $\S^0(t) \in \Z^N \cap \W_N, t \in \N_0$ is labeled.
We write the probability measure for $\Xi(t, \cdot), t \in \N_0$
started at $\xi \in \mM_0$ as
$\P_{\xi}$ with expectation $\E_{\xi}$,
and introduce a filtration
$\{\cF(t) : t \in \N_0\}$ defined by
$\cF(t)=\sigma(\Xi(s): 0 \leq s \leq t, s \in \N_0)$.
Then the above definition of the noncolliding RW gives the follows.  
Let $\xi=\sum_{j=1}^N \delta_{u_j}$ 
with $\u \in \Z^N_{\rm e} \cap \W_N$,
and $t \in \N$, $t \leq T \in \N$.
For any $\cF(t)$-measurable bounded function $F$,
\begin{equation}
\E_{\xi} \Big[F(\Xi(\cdot)) \Big]
= \lim_{n \to \infty} \rE_{\u} \left[ \left.
F \left(\sum_{j=1}^N \delta_{S_j(\cdot)} \right)
\right| \tau_{\u} > n \right].
\label{eqn:E1}
\end{equation}
The important fact is that, if we write
the Vandermonde determinant as
\begin{equation}
h(\x)=\det_{1 \leq j, k \leq N} [x_j^{k-1}]
=\prod_{1 \leq j < k \leq N} (x_k-x_j),
\label{eqn:Vand1}
\end{equation}
the expectation (\ref{eqn:E1}) is obtained by
an $h$-transform in the sense of Doob
of the form \cite{KOR02}
\begin{equation}
\E_{\xi}\Big[F(\Xi(\cdot))\Big]
=\rE_{\u} \left[ 
F \left(\sum_{j=1}^N \delta_{S_j(\cdot)} \right)
\1(\tau_{\u} > T)
\frac{h(\S(T))}{h(\u)} \right].
\label{eqn:E2}
\end{equation}
(See also \cite{Koe05,EK08}.) 
It determines the noncolliding RW,
$(\Xi(t), t \in \N_0, \P_{\xi})$.

The formula (\ref{eqn:E2}) is a discrete analogue
of the construction of noncolliding Brownian motion (BM)
by Grabiner \cite{Gra99} as an $h$-transform of
absorbing BM in $\W_N$.
The noncolliding BM is equivalent to Dyson's BM model
with parameter $\beta=2$ and the latter is known
as an eigenvalue process of Hermitian matrix-valued BM
and as solutions of the following system of
stochastic differential equations (SDEs)
\begin{equation}
dX_j(t)= dW_j(t)+
\sum_{\substack{1 \leq k \leq N, \cr k \not=j}}
\frac{1}{X_j(t)-X_k(t)}dt,
\quad 1 \leq j \leq N, \quad t \in [0, \infty),
\label{eqn:Dyson1}
\end{equation}
where $W_j(\cdot), 1 \leq j \leq N$
are independent one-dimensional standard BMs
\cite{Dys62,Meh04,Spo87,Gra99,Joh01,KT04,KT11b,Tao12,Osa12,GM13,Osa13,Osa13b}.
(From now on BM stands for one-dimensional standard 
Brownian motion and Dyson's BM model with $\beta=2$
is simply called {\it the Dyson model} in this paper.)
Then the noncolliding RW has been attracted much attention
as a discretization of models associated with
the Gaussian random-matrix ensembles
\cite{Bai00,Joh02,NF02,KT03,Joh05,BS07,For10,Fei12}.

Eigenvalue distributions of random-matrix ensembles
provide important examples of 
{\it determinantal point processes},
in which any correlation function is given by a determinant
specified by a single continuous function 
called the {\it correlation kernel}
\cite{Sos00,ST03,BKPV09}.
The noncolliding BM is regarded as a dynamical extension
of determinantal point process such that
any spatio-temporal correlation function
is expressed by a determinant.
Such processes are said to be {\it determinantal} \cite{KT07}.
The dynamical correlation kernel 
is asymmetric with respect to
the exchange of two points on the spatio-temporal plane
and shows causality in the system.
This type of correlation kernel 
was first obtained by Eynard and Mehta 
for a multi-matrix model \cite{EM98}
and by Nagao and Forrester 
for the noncolliding BM started at a special
initial distribution (the GUE eigenvalue distribution) \cite{NF98}.
It is proved that the noncolliding BM 
is determinantal for any fixed initial configuration with finite
numbers of particles as well as two families of 
infinite-particle initial configurations \cite{KT10,KT13}.

Nagao and Forrester \cite{NF02} studied
a `bridge' of noncolliding RW
started from $\u_0=(2j)_{j=0}^{N-1}$ at $t=0$
and returned to the same configuration $\u_0$
at time $t=2M, M \in \N_0$.
They showed that at time $t=M$
the spatial configuration provides 
a determinantal point process and the correlation 
kernel is expressed by using the symmetric Hahn polynomials.
Johansson \cite{Joh05} generalized the process to
a bridge from $\u_0$ at $t=0$
to $M_2-M_1+\u_0$ at $t=M_1+M_2$, $M_1, M_2 \in \N_0, M_2 > M_1$,
and proved that the process is determinantal.
The dynamical correlation kernel is of
the Eynard-Mehta type and called 
the extended Hahn-kernel.
For the noncolliding RW defined 
for infinite time-period $t \in \N_0$ by
(\ref{eqn:E1}) or (\ref{eqn:E2}) \cite{KOR02,Koe05,EK08},
however, determinantal structure of 
spatio-temporal correlations has not been clarified so far.

In the present paper we show that
the construction by the $h$-transform (\ref{eqn:E2})
directly leads to the fact that
the discrete-time noncolliding RW
is determinantal for any fixed initial configuration
$\xi=\sum_{j=1}^{N} \delta_{u_j} \in \mM_0$ with $N=\xi(\Z_{\rm e}) \in \N$.
(See \cite{Esa14} for the noncolliding system of 
continuous-time random walks.)
There are two key points in the present study of discrete-time systems;
proper {\it complexification} of RWs and introduction of 
{\it determinantal martingale}.
Let $\widetilde{\zeta} \in \R$ be a continuous random variable
in the {\it hyperbolic secant distribution} \cite{Fel66},
\begin{equation}
\widetilde{\rP}[\widetilde{\zeta} \in dx]
= \frac{1}{2} {\rm sech} \left( \frac{\pi x}{2} \right) dx
= \frac{1}{2 \cosh (\pi x/2)} dx,
\label{eqn:tP1}
\end{equation}
which is selfdecomposable (see pp.98-99 in \cite{Sat99}).
The Fourier transform of (\ref{eqn:tP1})
(the characteristic function of $\widetilde{\zeta}$) is also 
expressed by the hyperbolic secant \cite{Fel66}
($i \equiv \sqrt{-1}$)
\begin{equation}
\widetilde{\rE}[e^{i \alpha \widetilde{\zeta}}]
={\rm sech} \alpha
=\frac{1}{\cosh \alpha}, \quad \alpha \in \R,
\label{eqn:t_zeta2}
\end{equation}
which is exactly the inverse of (\ref{eqn:zeta1}).
Let $\{\widetilde{\zeta}(t) : t \in \N\}$
be a series of i.i.d.random variables obeying the same
probability law with $\widetilde{\zeta}$.
We define a discrete-time Markov process
$\widetilde{S}(t), t \in \N_0$ on $\R$
starting from 0 at time $t=0$ by
\begin{equation}
\widetilde{S}(t)=\widetilde{\zeta}(1) + \cdots + \widetilde{\zeta}(t),
\quad t \in \N.
\label{eqn:t_W1}
\end{equation}
At each time $t \in \N_0$, 
it is in the {\it generalized hyperbolic secant distribution}
with density
\begin{eqnarray}
\widetilde{p}(t, x|0)
&\equiv& \widetilde{\rP}[\widetilde{S}(t) \in dx ]
\nonumber\\
&=& \frac{2^{t-2}}{\pi \Gamma(t)}
\left| \Gamma \left( \frac{t}{2}+i \frac{x}{2} \right) \right|^2,
\quad t \in \N_0, \quad x \in \R,
\label{eqn:t_p1}
\end{eqnarray}
where $\Gamma$ denotes the gamma function
\cite{HH68}.
It can be shown that
$\widetilde{S}(t)/\sqrt{t} \dto {\sf N}(0,1)$ 
as $t \to \infty$ \cite{HH68}.
Let
$\widetilde{S}_j(\cdot), 1 \leq j \leq N$
be a set of independent 
copies of $\widetilde{S}(\cdot)$ and 
express the expectation with respect to these processes
also by $\widetilde{\rE}$.
For the original RW, 
$\S(t)=(S_1(t), \dots, S_N(t)), t \in \N_0$
started at a fixed configuration 
$\u \in \Z_{\rm e}^N \cap \W_N$,
its {\it complexification} is given by 
the discrete-time complex processes,
$\bZ(t)=(Z_1(t), \dots, Z_N(t)), t \in\N_0$ with
\begin{equation}
Z_j(t)=S_j(t)+ i \widetilde{S}_j(t),
\quad 1 \leq j \leq N, \quad t \in \N_0.
\label{eqn:Z1}
\end{equation}
We put $\xi = \sum_{j=1}^{N} \delta_{u_j} \in \mM_0$
and consider a set of functions of $z \in \C$, 
\begin{eqnarray}
\Phi_{\xi}^{u_k}(z)
&=& \prod_{\substack{1 \leq j \leq N,  \cr j \not= k}}
\frac{z-u_j}{u_k-u_j},
\quad 1 \leq k \leq N.
\label{eqn:Phi1}
\end{eqnarray}
The function $\Phi_{\xi}^{u_k}(z)$ is a polynomial
of $z$ with degree $N-1$ having zeros
at $u_j, 1 \leq j \leq N, j \not=k$
and $\Phi_{\xi}^{u_k}(u_k)=1$.
We can prove that (Lemma \ref{thm:poly_mart}), 
for each $1 \leq k \leq N$,
\begin{equation}
\cM_{\xi}^{u_k}(t, S_j(t))
\equiv \widetilde{\rE}
\Big[ \Phi_{\xi}^{u_k}(Z_j(t)) \Big],
\quad 1 \leq j \leq N
\label{eqn:mart1}
\end{equation}
provide independent martingales with discrete time $t \in \N_0$.
We consider a determinant of matrix, whose entries
are these martingales,
\begin{equation}
\cD_{\xi}(t,\S(t))
=\det_{1 \leq j, k \leq N}
[\cM_{\xi}^{u_k}(t, S_j(t))],
\quad t \in \N_0,
\label{eqn:mart2}
\end{equation}
which we call the
{\it determinantal martingale} \cite{Kat14}.
Our martingales (\ref{eqn:mart1}) are prepared so that 
the equality 
\begin{equation}
\frac{h(\S(t))}{h(\u)}
=\cD_{\xi}(t, \S(t)), \quad t \in \N_0, 
\label{eqn:dm2}
\end{equation}
holds and a kind of reducibility (Lemma \ref{thm:reducibility}) is
established. 

This equality (\ref{eqn:dm2}) gives
a {\it determinantal-martingale representation (DMR)}
for the noncolliding RW (Proposition \ref{thm:DMR1}), 
and from it we can prove that the noncolliding RW is determinantal
with the correlation kernel,
\begin{equation}
\mbK_{\xi}(s,x;t,y)
= \left\{ \begin{array}{l}
\displaystyle{
\sum_{j=1}^N p(s, x|u_j) \cM_{\xi}^{u_j}(t,y)
- \1(s>t) p(s-t,x|y),
}
\cr
\qquad \qquad \qquad \mbox{if} \quad
(s,x), (t,y) \in \N_0 \times \Z, \quad s+x, t+y \in \Z_{\rm e},
\cr
0,
\qquad \qquad \quad \mbox{otherwise},
\end{array} \right.
\label{eqn:K1}
\end{equation}
where $p$ is the transition probability (\ref{eqn:tp1}), 
and 
$\1(\cdot)$ is an indicator;
$\1(\omega)=1$ if $\omega$ is satisfied,
and $\1(\omega)=0$ otherwise (Theorem \ref{thm:main1}).
Note again that 
\begin{eqnarray}
\cM_{\xi}^{u_j}(t,y)
&=& \widetilde{\rE} \left[ 
\Phi_{\xi}^{u_j}(y+i \widetilde{S}(t)) \right]
\nonumber\\
&=& \int_{\R} dv \, \widetilde{p}(t, v|0) 
\Phi_{\xi}^{u_j}(y+i v), \quad 1 \leq j \leq N
\label{eqn:mart1b}
\end{eqnarray}
with (\ref{eqn:t_p1}), 
are functions of initial configuration
$\xi=\sum_{j=1}^N \delta_{u_j}$ through (\ref{eqn:Phi1}).

For $a \in \{2,3, \dots\}$, we consider a configuration on $\Z_{\rm e}$ having
equidistant spacing $2a$ with an infinite number of particles,
\begin{equation}
\delta_{2a \Z}(\cdot) \equiv \sum_{k \in \Z} \delta_{2a k} (\cdot).
\label{eqn:xi2a}
\end{equation}
(The noncolliding RW starting from
$\delta_{2 \Z}(\cdot)$, that is,
the case $a=1$ of (\ref{eqn:xi2a}),
is trivial.
The process is stationary in the sense that
$\Xi(2n)= \sum_{k \in \Z} \delta_{2k}$,
$\Xi(2n+1)= \sum_{k \in \Z} \delta_{2k+1}$,
$n \in \N_0$.)
We prove that the noncolliding RW started at (\ref{eqn:xi2a}), 
denoted as $(\Xi(t), t \in \N_0, \P_{\delta_{2a \Z}})$,
$a \in \{2,3, \dots \}$,
is well-defined as a determinantal process
with an infinite number of particles (Proposition \ref{thm:inf_det}).
There the $N$ linearly independent 
polynomials of $y$ given by (\ref{eqn:mart1b})
are extended to an infinite sequence of 
linearly independent entire functions of $y$,
$\cM_{\delta_{2a \Z}}^{2ak}(t,y), k \in \Z$,
corresponding to the infinite-particle 
initial configuration (\ref{eqn:xi2a}).
Then by using the infinite sequence of independent
martingales with discrete time, 
$(\cM_{\delta_{2a \Z}}^{2ak}(t, S_j(t)))_{t \in \N_0}, k \in \Z$,
for each $j \in \Z$, 
we can give DMRs for $(\Xi(t), t \in \N_0, \P_{\delta_{2 a \Z}})$,
$a \in \{2,3, \dots\}$ (Proposition \ref{thm:inf_noncRW}).
For each $a \in \{2,3, \dots\}$, 
this discrete-time infinite-particle system on $\Z$ 
shows a relaxation phenomenon to the equilibrium
determinantal process, $(\Xi(t), t \in \Z, \bP_{\rho}$), 
whose correlation kernel is given by
\begin{equation}
\bK_{\rho}(t-s, y-x) = \left\{
\begin{array}{ll}
\displaystyle{\int_0^{\rho} d u \,
\frac{2 \cos(\pi u(y-x))}{[\cos(\pi u)]^{t-s}}},
\quad & \mbox{if $s<t$}, \cr
& \cr
\displaystyle{\frac{2 \sin(\pi \rho(y-x))}{\pi(y-x)}},
\quad & \mbox{if $s=t$}, \cr
& \cr
\displaystyle{-\int_{\rho}^{1} d u \,
\frac{2 \cos(\pi u(y-x))}{[\cos(\pi u)]^{t-s}}},
\quad & \mbox{if $s>t$},
\end{array}
\right.
\label{eqn:K3b}
\end{equation}
for $(s,x), (t,y) \in \Z^2$,
$s+x, t+y \in \Z_{\rm e}$,
and $\bK_{\rho}(t-s, y-x)=0$ otherwise, 
where $\rho=1/2a$ gives the particle density on $\Z$
(Theorem \ref{thm:equilibrium}).
This is a discrete analogue of the
extended sine-kernel (see Section 11.7.1 in \cite{For10}) of 
the Dyson model (\ref{eqn:Dyson1}). 

We note that independent increments
$\zeta_j(t)$ of $S_j(t)$ and
$\widetilde{\zeta}_j(t)$ of $\widetilde{S}_j(t)$,
$1 \leq j \leq N, t \in \N_0$
are both having mean zero and variance 1.
Then Donsker's {\it invariance principle}
\cite{Bil68,RY05} proves
both of $S_j(n^2 t)/n$ and $\widetilde{S}_j(n^2 t)/n$
converge to BMs as $n \to \infty$.
It implies that the DMRs for appropriately scaled
noncolliding RWs converge to 
the {\it complex BM representation} for the
Dyson model (\ref{eqn:Dyson1}) given by \cite{KT13}.
The central limit theorem of noncolliding RWs
to the Dyson model will be established.

The paper is organized as follows.
In Section \ref{sec:Preliminaries}
the polynomial martingales and determinantal martingales
are introduced for noncolliding RW 
and their properties are discussed.
Determinantal properties of noncolliding RW 
is clarified in Section \ref{sec:det_prop}.
An extension to infinite particle systems is discussed in 
Section \ref{sec:infinite}.
Convergence of noncolliding RWs
to the Dyson model 
is discussed in Section \ref{sec:convergence}.

\SSC{Preliminaries \label{sec:Preliminaries}}
\subsection{Discrete It\^o's formula 
and polynomial martingales of Fujita}

Let $S(t), t \in \N_0$ be a one-dimensional, 
simple and symmetric RW 
starting from 0 at time $t=0$,
$$
S(t)=\zeta(1)+\zeta(2)+ \cdots + \zeta(t),
\quad t \in \N,
$$
where $\{\zeta(t) : t \in \N\}$ are i.i.d.obeying 
the same probability law with $\zeta$.
The following {\it discrete It\^o's formula}
was given by Fujita 
for the one-dimensional, simple and symmetric RW \cite{Fuj02,Fuj08}.

\begin{lem}
\label{thm:Fujita1}
For any $f: \N_0 \times \Z \to \R$ and
any $t \in \N_0$,
\begin{eqnarray}
&& f(t+1,S(t+1))-f(t,S(t))
\nonumber\\
&& \quad = \frac{1}{2} \Big[
f(t+1, S(t)+1)-f(t+1, S(t)-1) \Big] \zeta(t+1)
\nonumber\\
&& \qquad + \frac{1}{2} \Big[ f(t+1, S(t)+1)- 2 f(t+1, S(t))+ f(t+1, S(t)-1) \Big]
\nonumber\\
&& \qquad +f(t+1, S(t))-f(t, S(t)).
\label{eqn:Ito_Fujita}
\end{eqnarray}
\end{lem}
\vskip 0.3cm

We perform the Esscher transform 
with parameter $\alpha \in \R$,
$S(\cdot) \to \widehat{S}_{\alpha}(\cdot)$ as
\begin{equation}
\widehat{S}_{\alpha}(t)
=\frac{e^{\alpha S(t)}}{\rE[e^{\alpha S(t)}]},
\quad t \in \N_0.
\label{eqn:tV1}
\end{equation}
By (\ref{eqn:zeta1}), 
we have
$
\widehat{S}_{\alpha}(t)=
G_{\alpha}(t, S(t))
$
with
\begin{equation}
G_{\alpha}(t,x)=\frac{e^{\alpha x}}
{(\cosh \alpha)^t},
\quad t \in \N_0, \quad x \in \Z.
\label{eqn:G1}
\end{equation}
If we set $f=G_{\alpha}$ in (\ref{eqn:Ito_Fujita}),
the second and third terms in the RHS vanish.
Then
\begin{eqnarray}
&& G_{\alpha}(t+1, S(t+1))-G_{\alpha}(t,S(t))
\nonumber\\
&& \qquad
= \frac{1}{2} \Big[
G_{\alpha}(t+1, S(t)+1) - G_{\alpha}(t+1, S(t)-1) \Big]
\zeta(t+1),
\label{eqn:G2}
\nonumber
\end{eqnarray}
which implies that $G_{\alpha}(t, S(t))$ is
$\{\zeta(1), \dots, \zeta(t)\}$-martingale
for any $\alpha \in \R$ \cite{Fuj02,Fuj08}.
From now on, we simply say
`$(G_{\alpha}(t, S(t)))_{t \in \N_0}$ is a martingale'
in such a situation.

Expansion of (\ref{eqn:G1}) with respect to $\alpha$
around $\alpha=0$, 
\begin{equation}
G_{\alpha}(t,x)=\sum_{n=0}^{\infty} m_n(t,x)
\frac{\alpha^n}{n!},
\label{eqn:G3}
\end{equation}
determines a series of monic polynomials
of degrees $n$ studied by Fujita in \cite{Fuj02}
\begin{equation}
m_n(t,x)=x^n + \sum_{j=1}^{n-1} c_n^{(j)}(t) x^j,
\quad n \in \N_0, 
\label{eqn:m1}
\end{equation}
such that
\begin{eqnarray}
\label{eqn:cond1}
&& c_n^{(j)}(0)=0, \quad 1 \leq j \leq n-1, 
\\
\label{eqn:cond2}
&&
\mbox{$m_n(t, S(t))$ is martingale}, \quad
t \in \N_0.
\end{eqnarray}
For example,
\begin{eqnarray}
m_0(t,x) &=& 1,
\nonumber\\
m_1(t,x) &=& x,
\nonumber\\
m_2(t,x) &=& x^2-t,
\nonumber\\
m_3(t,x) &=& x^3-3 t x,
\nonumber\\
m_4(t,x) &=& x^4-6 t x^2+t (3t+2),
\nonumber\\
m_5(t,x) &=& x^5-10 t x^3 + 5 t (3t+2) x,
\quad \cdots .
\label{eqn:mn_ex}
\end{eqnarray}
They satisfy the recurrence relations
$$
m_n(t,x)=\frac{1}{2}
[ m_n(t+1,x+1)+m_n(t+1,x-1) ],
\quad n \in \N_0.
$$
As mentioned below in Remark 2,
$m_n(t, x), n \in \N_0$ are related with the Euler polynomials
studied in \cite{Sch00}.
Since the importance of $m_n(t, x), n \in \N_0$ 
in the context of random walks was first clearly shown 
by Fujita \cite{Fuj02}, however, 
we would like to call $m_n(t,x), n \in \N$, 
{\it Fujita's polynomials} and 
$(m_n(t,S(t)))_{t \in \N_0}, n \in \N_0$, 
{\it Fujita's polynomial martingales} 
for the simple and symmetric RW.

\vskip 0.3cm
\noindent{\bf Remark 1.} \,
Let $B(t), t \geq 0$ be BM started at 0.
Then its Esscher transform with parameter $\alpha$ is given by
$
\widehat{B}_{\alpha}(t)
= G^{\rm BM}_{\alpha}(t,B(t))
$
with
$$
G^{\rm BM}_{\alpha}(t,x)
= \frac{e^{\alpha x}}{\rE[e^{\alpha B(t)}]}
=\frac{e^{\alpha x}}
{\displaystyle{ \int_{-\infty}^{\infty} dx
e^{\alpha x} p^{\rm BM}(t, x|0)
 }}
= e^{\alpha x - \alpha^2 t/2},
$$
where 
\begin{equation}
p^{\rm BM}(t,y|x)
=\frac{1}{\sqrt{2 \pi t}} e^{-(y-x)^2/2t},
\quad
t \geq 0, \quad x, y \in \R
\label{eqn:pBM}
\end{equation}
is the transition probability density of BM.
We see that
$$
G^{\rm BM}_{\alpha}(t,x)
=\sum_{n=0}^{\infty} \left( \frac{t}{2} \right)^{n/2}
H_n \left( \frac{x}{\sqrt{2t}}\right) \frac{\alpha^n}{n!}
$$
with the Hermite polynomials
$\displaystyle{
H_n(z)= \sum_{j=0}^{[n/2]} (-1)^j
\frac{n !}{j ! (n-2j)!} (2z)^{n-2j},
n \in \N_0.
}$
Therefore, 
\begin{equation}
m_n^{\rm BM}(t, B(t))
=\left( \frac{t}{2} \right)^{n/2}
H_n \left( \frac{B(t)}{\sqrt{2t}} \right),
\quad n \in \N_0, \quad t \geq 0,
\label{eqn:mnBM1}
\end{equation}
are the polynomial martingales for BM 
as known well (see, for instance, \cite{Sch00}).
\vskip 0.3cm

\vskip 0.3cm
\noindent{\bf Remark 2.} \,
The polynomials (\ref{eqn:mnBM1}) for BM have
the multiple stochastic-integral representations,
$$
m_n^{\rm BM}(t, B(t))
= n! \int_0^t \int_0^{t_1} \cdots \int_0^{t_{n-1}} 
dB(t_n) \cdots d B(t_2) d B(t_1),
\quad n \in \N.
$$
Their discrete analogues determine the polynomial martingales
for RW,
$$
\widehat{m}_n(t, S(t))
= n! \sum_{t_1=1}^{t} \sum_{t_2=1}^{t_1} \cdots
\sum_{t_n=1}^{t_{n-1}}
\zeta(t_1) \zeta(t_2) \cdots \zeta(t_n).
$$
For $N \in \N_0$, $0 < p < 1$,  
the monic Krawtchouk polynomials 
$\tilde{K}_n(x;N,p), n \in \N_0$ are defined by the generating function as
$$
\sum_{n=0}^{N} \tilde{K}_n(x;N,p) \frac{\alpha^n}{n!}
=(1+(1-p) \alpha)^{x}(1-p \alpha)^{N-x}.
$$
Then \cite{Sch00,Pri09}
$$
\widehat{m}_n(t,x)=2^n \tilde{K}_n((t+x)/2; t, 1/2),
\quad n \in \N_0.
$$
It should be noted that
$\widehat{m}_n(t,x), n \in \N_0$ are generally
different from Fujita's polynomials $m_n(t,x), n \in \N_0$.
We see that
$\widehat{m}_0(t,x)=1$,
$\widehat{m}_1(t,x)=x$,
$\widehat{m}_2(t,x)=x^2-t$, and
$\widehat{m}_3(t,x)=x(x^2+2)-3tx$,
$\widehat{m}_4(t,x)=x^2(x^2+8)-6tx^2+3t(t-2),
\dots$.
In general, the Krawtchouk polynomials do not satisfy 
the condition (\ref{eqn:cond1}).
The monic polynomials of order $n$, $E_n^{(\lambda)}(x)$ with 
parameter $\lambda \in \N_0$ defined by the generating function
\begin{equation}
\sum_{n=0}^{\infty} E_n^{(\lambda)}(x) \frac{\alpha^n}{n!}
=\left( \frac{2}{1+e^{\alpha}} \right)^{\lambda} e^{\alpha x}
\label{eqn:Euler1}
\end{equation}
are called the Euler polynomials 
(see p.253 in \cite{Erd55}).
Schoutens showed that, if $\bar{\zeta}_j, j \in \N$ have 
a binomial distribution ${\sf Bin}(\lambda, 1/2)$ 
and $\bar{S}(t) \equiv \sum_{j=1}^t \bar{\zeta}_j, t \in \N$ with $\bar{S}(0) \equiv 0$,
then $(E_n^{(t \lambda)}(\bar{S}(t)))_{t \in \N_0}, n \in \N_0$
are martingales \cite{Sch00}.
Fujita's polynomials are related with Euler's by
\begin{equation}
m_n(t, x)=2^n E_n^{(t)} \left( \frac{t+x}{2} \right), 
\quad n \in \N_0.
\label{eqn:Euler2}
\end{equation}

\subsection{Complex-process representation for polynomial martingales}
\label{sec:CPR_poly}

For RW, $(S(t))_{t \in \N_0}$, we consider 
its complexification, 
\begin{equation}
Z(t)=S(t)+i \widetilde{S}(t), \quad t \in \N_0,
\label{eqn:cp2}
\end{equation}
where $\widetilde{S}(\cdot)$ is defined
by (\ref{eqn:t_W1}) with $\widetilde{S}(0) \equiv 0$.
Note that $\Re Z(t)=S(t) \in \Z$ and
$\Im Z(t)=\widetilde{S}(t) \in \R, t \in \N_0$.
We can prove the following.

\begin{lem}
\label{thm:cpr}
With the complex process (\ref{eqn:cp2}), 
Fujita's polynomial martingales with discrete time $t \in \N_0$, 
$(m_n(t,S(t)))_{t \in \N_0}$, $n \in \N_0$, for the
simple and symmetric RW have
the following representations,
\begin{equation}
m_n(t, S(t))
= \widetilde{\rE} [ Z(t)^n ],
\quad n \in \N_0, \quad t \in \N_0.
\label{eqn:cp3}
\end{equation}
\end{lem}
\noindent{\it Proof.} \,
By definition (\ref{eqn:t_W1}) of $\widetilde{S}(t)$,
(\ref{eqn:t_zeta2}) gives
\begin{equation}
\widetilde{\rE} \left[ e^{i \alpha \widetilde{S}(t)} \right]
=\left( \widetilde{\rE} \left[e^{i \alpha \widetilde{\zeta}} \right] \right)^t
=\frac{1}{(\cosh \alpha)^t},
\quad \alpha \in \R, \quad t \in N_0.
\label{eqn:t_W2}
\end{equation}
Then for (\ref{eqn:G1}), the equality
$
G_{\alpha}(t, S(t)) = 
\widetilde{\rE}[e^{\alpha Z(t)}],
\alpha \in \R
$
is established, 
which proves (\ref{eqn:cp3}). \qed
\vskip 0.3cm

\vskip 0.3cm
\noindent{\bf Remark 3.} \,
For a pair of independent BMs,
$B(t), \widetilde{B}(t), t \geq 0$,
we can see
\begin{equation}
\rE[e^{\alpha B(t)}]
= e^{\alpha^2 t/2}
=\left( \widetilde{\rE}
\left[ e^{i \alpha \widetilde{B}(t)} \right] \right)^{-1},
\quad \alpha \in \R.
\label{eqn:recip1}
\end{equation}
Then
$
m_n^{\rm BM}(t, B(t))=\widetilde{\rE} [\cB(t)^n], 
n \in \N_0,t \geq 0,
$
is concluded, 
where $\cB(t)$ is a complex BM,
$\cB(t)=B(t)+i \widetilde{B}(t), t \geq 0$.
The reciprocity relations between
(\ref{eqn:zeta1}) and (\ref{eqn:t_zeta2}),
and $\rE[e^{\alpha S(t)}]=(\cosh \alpha)^t$
and (\ref{eqn:t_W2})
are discrete-time analogues of (\ref{eqn:recip1}).
\vskip 0.3cm

A direct consequence of Lemma \ref{thm:cpr} is the following.
\begin{lem}
\label{thm:poly_mart}
Assume that $f$ is polynomial. Then
$\widetilde{\rE}[f(Z(t))]$ is a martingale with discrete time $t \in \N_0$.
\end{lem}

\subsection{Determinantal martingales}

We consider an $N$-component complex process
$\bZ(t)=(Z_1(t), \dots, Z_N(t)), t \in \N_0$ with (\ref{eqn:Z1}).
The probability space for (\ref{eqn:Z1}) is a product
of the probability space $(\Omega, \cF, \rP_{\u})$ 
for the RW on $\Z^N$, $\S(t), t \in \N_0$, 
and $(\widetilde{\Omega}, \widetilde{\cF}, \widetilde{\rP})$
for $\widetilde{\S}(t), t \in \N_0$ defined on $\R^N$,
which we write as $(\check{\Omega}, \check{\cF}, \bP_{\u})$.
Let $\bE_{\u}$ be the expectation for the process $\bZ(t), t \in \N_0$ 
with the initial condition $\bZ(0)=\u \in \Z_{\rm e}^N \cap \W_N$.

By multilinearity of determinant, the Vandermonde determinant 
(\ref{eqn:Vand1}) does not change
in replacing $x_j^{k-1}$ by any monic polynomial of $x_j$ of 
degree $k-1$, $1 \leq j, k \leq N$.
Note that $m_{k-1}(t,x_j)$ is a monic polynomial of $x_j$ of degree $k-1$. 
Then
\begin{eqnarray}
\frac{h(\S(t))}{h(\u)}
&=& \frac{1}{h(\u)} \det_{1 \leq j, k \leq N} [m_{k-1}(t, S_j(t))]
\nonumber\\
&=& \frac{1}{h(\u)} \det_{1 \leq j, k \leq N}
\Big[ \widetilde{\rE} [Z_j(t)^{k-1}  ] \Big]
\nonumber\\
&=& \widetilde{\rE} \left[ \frac{1}{h(\u)}
\det_{1 \leq j, k \leq N}
[Z_j(t)^{k-1}] \right],
\nonumber
\end{eqnarray}
where we have used Lemma \ref{thm:cpr}, 
the multilinearity of determinant, and
independence of $Z_j(t)$'s.
Therefore, we have obtained the equality,
\begin{equation}
\frac{h(\S(t))}{h(\u)}
= \widetilde{\rE} \left[
\frac{h(\bZ(t))}{h(\u)} \right], \quad t \in \N_0.
\label{eqn:Vand2}
\end{equation}

Now we consider the determinant identity \cite{KT13},
\begin{equation}
\frac{h(\z)}{h(\u)}
=\det_{1 \leq j, k \leq N} \Big[
\Phi_{\xi}^{u_k}(z_j) \Big],
\label{eqn:Vand3}
\end{equation}
where $\xi=\sum_{j=1}^N \delta_{u_j}, \u =(u_1, \dots, u_N) \in \W_N$, and
$\Phi_{\xi}^{u_k}(z)$ is given by
(\ref{eqn:Phi1})
(see Section 4.1 in \cite{Kat14} for derivation). 
Since $\Phi_{\xi}^{u_k}(z)$ is a polynomial of $z$
of degree $N-1$, 
Lemma \ref{thm:poly_mart} proves that 
$(\cM_{\xi}^{u_k}(t, S_j(t)))_{t \in \N_0}, 1 \leq j \leq N$,
defined by (\ref{eqn:mart1})
are independent martingales with discrete time $t \in \N_0$ and 
\begin{eqnarray}
\rE_{\u}[\cM_{\xi}^{u_k}(t, S_j(t))]
&=& \rE_{\u}[\cM_{\xi}^{u_k}(0, S_j(0))]
\nonumber\\
&=& \cM_{\xi}^{u_k}(0, u_j)
\nonumber\\
&=& \Phi_{\xi}^{u_k}(u_j)=\delta_{jk},
\quad 1 \leq j, k \leq N.
\label{eqn:cM2}
\end{eqnarray} 
Using the identity (\ref{eqn:Vand3})
for $h(\bZ(t))/h(\u)$ in (\ref{eqn:Vand2}), 
we have
\begin{eqnarray}
\frac{h(\S(t))}{h(\u)}
&=& \widetilde{\rE} \left[
\det_{1 \leq j, k \leq N}
[ \Phi_{\xi}^{u_k}(Z_j(t)) ] \right]
\nonumber\\
&=& \det_{1 \leq j, k \leq N} \Big[
\widetilde{\rE} [ \Phi_{\xi}^{u_k}(Z_j(t)) ] \Big],
\nonumber
\end{eqnarray}
where independence of $Z_j(t)$'s was again used.
By definition (\ref{eqn:mart2}) of $\cD_{\xi}$ with (\ref{eqn:mart1}), 
we obtain the equality (\ref{eqn:dm2}).

\vskip 0.3cm
\noindent{\bf Remark 4.} \,
The real parts of the complex processes
(\ref{eqn:Z1})
are RWs with 
$\rE_{\u}[(S_j(t)-u_j)^2]=t \in \N_0, 1 \leq j \leq N$.
It is obvious from definition (\ref{eqn:t_W1}) that
the imaginary parts, $\widetilde{S}_j(t), t \in \N_0$, 
are $\{\widetilde{\zeta}_j(1), \dots, \widetilde{\zeta}_j(t) \}$-martingales
with $\widetilde{\rE}[\widetilde{S}_j(t)^2]=t \in \N_0,1 \leq j \leq N$.
Then $Z_j(\cdot), 1 \leq j \leq N$ 
shall be regarded as {\it discrete-time conformal martingales}
(see Definition (2.2) in Section V.2 of \cite{RY05}).
Their conformal maps by polynomial functions,
$\Phi_{\xi}^{u_k}(Z_j(\cdot)), 1 \leq j, k \leq N$ are 
discrete-time complex martingales such that
\begin{equation}
\bE_{\u} [\Phi_{\xi}^{u_k}(Z_j(t))]
=\bE_{\u}[\Phi_{\xi}^{u_k}(Z_j(0))] =\delta_{j k},
\quad 1 \leq j, k \leq N
\label{eqn:conformal1}
\end{equation}
for any $t \in \N_0$, 
\vskip 0.3cm

For $n \in \N$, let
$\I_{n} = \{1,2, \dots, n\}$.
Denote the cardinality of a finite set $A$
by $|A|$.
Let  $\x=(x_1, \dots, x_N) \in \Z^N$
and $1 \leq N' \leq N$.
We write $\J \subset \I_N, |\J|=N'$,
if $\J=\{j_1, \dots, j_{N'}\},
1 \leq j_1 < \dots < j_{N'} \leq N$,
and put $\x_{\J}=(x_{j_1}, \dots, x_{j_{N'}})$.
In particular, we write
$\x_{N'}=\x_{\I_{N'}}, 1 \leq N' \leq N$.
(By definition $\x_N=\x$.)
Suppose $\u \in \Z_{\rm e}^N \cap \W_N$ and 
$\xi(\cdot)=\sum_{j=1}^N \delta_{u_j}(\cdot)$.
For $\J \subset \I_N,  |\J|=N', 1 \leq N' \leq N$,
introduce determinantal martingales
\begin{equation}
\cD_{\xi}(t, \S_{\J}(t)) 
=\det_{j,k \in \J} \Big[ 
\cM_{\xi}^{u_k}(t, S_j(t)) \Big],
\quad t \in \N_0,
\label{eqn:dmB1}
\end{equation}
where the sizes of matrices for determinants 
are $|\J|=N', 1 \leq N' \leq N$.
We can prove the following,
which is a discrete-time version of Lemma 2.1 in \cite{Kat14}.

\begin{lem}
\label{thm:reducibility}
Assume that
$\xi(\cdot)=\sum_{j=1}^N \delta_{u_j}(\cdot)$
with $\u \in \Z_{\rm e}^N \cap \W_N$.
Let $1 \leq N' < N$.
For $t \in \N, t \leq T \in \N$ and 
a symmetric bounded function
$F_{N'}$ on $\Z^{N'}$,
\begin{eqnarray}
&& \sum_{\J \subset \I_N, |\J|=N'}
\rE_{\u} \left[
F_{N'}(\S_{\J}(t))
\cD_{\xi}(T, \S(t)) \right]
\nonumber\\
&& \quad
= \int_{\W_{N'}} \xi^{\otimes N'} (d\v)
\rE_{\v} \left[
F_{N'}(\S_{N'}(t))
\cD_{\xi}(T, \S_{N'}(T)) \right].
\label{eqn:reducibility}
\end{eqnarray}
\end{lem}
This shows the
{\it reducibility} of the determinantal martingale
in the sense that,
if we observe a symmetric function $F_{N'}$ depending
on only $N'$ variables, $N' < N$,
then the size of determinant for determinantal
martingale can be reduced from $N$ to $N'$.

\SSC{Determinantal Properties \label{sec:det_prop}}
\subsection{Determinantal martingale representation}

Since we consider the noncolliding RW as a process
represented by an unlabeled configuration (\ref{eqn:Xi1}),
measurable functions of $\Xi(\cdot)$ are only
symmetric functions of $N$ variables, 
$S^0_j(\cdot), 1 \leq j \leq N$.
Then by the equality (\ref{eqn:dm2}), we obtain the
following representation.
We call it the {\it determinantal-martingale representation} (DMR)
for the present noncolliding RW.

\begin{prop}
\label{thm:DMR1}
Suppose that $N \in \N$ and
$\xi=\sum_{j=1}^{N} \delta_{u_j}$
with $\u=(u_1, \dots, u_N) \in \Z_{\rm e}^N \cap \W_N$.
Let $t  \in \N, t \leq T \in \N$.
For any ${\cal F}(t)$-measurable bounded function $F$
we have
\begin{eqnarray}
\E_{\xi} \left[ F \left(\Xi(\cdot) \right) \right]
&=& \rE_{\u} \left[F \left( \sum_{j=1}^{N} \delta_{S_j(\cdot)} \right)
\cD_{\xi}( T, \S(T)) \right]
\nonumber\\
&=& \bE_{\u} \left[F \left( \sum_{j=1}^{N} \delta_{\Re Z_j(\cdot)} \right)
\det_{1 \leq j, k \leq N}
[\Phi_{\xi}^{u_k}(Z_j(T))] \right]. 
\label{eqn:DMR1}
\end{eqnarray}
\end{prop}
\vskip 0.3cm
Note that the second representation of (\ref{eqn:DMR1})
is a discrete-time analogue of the {\it complex BM representation}
reported in \cite{KT13} for the Dyson model
({\it i.e.} the noncolliding BM).
See Remark 4 above again.
\vskip 0.3cm
\noindent{\it Proof of Proposition \ref{thm:DMR1}.} \,
It is sufficient 
to consider the case that $F$ is given as
$F(\Xi(\cdot))
= \prod_{m=1}^M g_m(\S^0(t_m))$
for $M \in \N$, $t_m \in \N, 1 \leq m \leq M$, 
$t_1< \cdots <t_M \leq T \in \N$, 
with symmetric bounded functions $g_m$ on 
$\Z^N$, $1 \leq m \leq M$.
Here we prove the equalities
\begin{eqnarray}
\E_{\xi}\left[ \prod_{m=1}^M g_m(\S^0(t_m)) \right]
&=& \rE_{\u} \left[ \prod_{m=1}^M g_m(\S(t_m))
\cD_{\xi}( T, \S(T)) \right]
\nonumber\\
&=& \bE_{\u} \left[ \prod_{m=1}^M g_m(\S(t_m))
\det_{1 \leq j, k \leq N}
[\Phi_{\xi}^{u_k}(Z_j(T))] \right].
\label{eqn:prA1}
\end{eqnarray}
By (\ref{eqn:E2}), the LHS of (\ref{eqn:prA1}) is given by
\begin{equation}
\rE_{\u} \left[ 
\prod_{m=1}^M g_m(\S(t_m))
\1(\tau_{\u} > t_M)
\frac{h(\S(t_M))}{h(\u)} \right],
\label{eqn:prA2}
\end{equation}
where we used the fact that $h(\S(\cdot))/h(\u)$ is martingale.
At time $t=\tau_{\u}$, there are at least one pair $(j, j+1)$ 
such that $S_j(\tau_{\u})=S_{j+1}(\tau_{\u}), 1 \leq j \leq N-1$.
We choose the minimal $j$.
Let $\sigma_{j, j+1}$ be the permutation of
the indices $j$ and $j+1$ and for
$\v=(v_1, \dots, v_N) \in \Z^N$ we put
$\sigma_{j,j+1}(\v)=(v_{\sigma_{j,j+1}(k)})_{k=1}^N
=(v_1, \dots, v_{j+1}, v_j, \dots, v_N)$.
Let $\u'$ be the labeled configuration of the process
at time $t=\tau_{\u}$.
Since $u'_j=u'_{j+1}$ by the above setting,
under the probability law $\rP_{\u'}$
the processes $\S(t), t > \tau_{\u}$ and 
$\sigma_{j,j+1}(\S(t)), t > \tau_{\u}$
are identical in distribution.
Since $g_m, 1 \leq m \leq M$ are symmetric,
but $h$ is antisymmetric, the Markov property
of the process $\S(\cdot)$ gives
$$
\rE_{\u} 
\left[ 
\prod_{m=1}^M g_m(\S(t_m))
\1(\tau_{\u} \leq t_M)
\frac{h(\S(t_M))}{h(\u)} \right]=0.
$$
Therefore, (\ref{eqn:prA2}) is equal to
$$
\rE_{\u} \left[ 
\prod_{m=1}^M g_m(\S(t_m))
\frac{h(\S(t_M))}{h(\u)} \right].
$$
By the equality (\ref{eqn:dm2}) and the
martingale property of $(\cD_{\xi}(t, \S(t)))_{t \in \N_0}$, 
we obtain the first line of (\ref{eqn:prA1}). 
By definitions of $\bE_{\u}$ and $\cD_{\xi}$, the second line is valid.
Then the proof is completed.
\qed
\vskip 0.5cm

\subsection{Determinantal process \label{sec:det_pr}}

For any integer $M \in \N$,
a sequence of times
$\t=(t_1,\dots,t_M) \in \N^M$ with 
$t_1 < \cdots < t_M \leq T \in \N$,
and a sequence of bounded functions
$\f=(f_{t_1},\dots,f_{t_M})$,
the {\it moment generating function} of multitime distribution
of the process $\Xi(\cdot)$ is defined by
\begin{equation}
\Psi_{\xi}^{\t}[\f]
\equiv \E_{\xi} \left[ \exp \left\{ \sum_{m=1}^{M} 
\int_{\Z} f_{t_m}(x) \Xi(t_m, dx) \right\} \right].
\label{eqn:GF1}
\end{equation}
It is expanded with respect to
\begin{equation}
\chi_{t_m}(\cdot)=e^{f_{t_m}(\cdot)}-1,
\quad
1 \leq m \leq M
\label{eqn:chi1}
\end{equation}
as
\begin{equation}
\Psi_{\xi}^{\t}[\f]
=\sum_
{\substack
{N_m \geq 0, \\ 1 \leq m \leq M} }
\sum_
{\substack
{\x^{(m)}_{N_m} \in \Z^{N_m} \cap \W_{N_m}, 
\\ 1 \leq m \leq M} }
\prod_{m=1}^{M}
\prod_{j=1}^{N_{m}} 
\chi_{t_m} \Big(x_{j}^{(m)} \Big)
\rho_{\xi} 
\Big( t_{1}, \x^{(1)}_{N_1}; \dots ; t_{M}, \x^{(M)}_{N_M} \Big),
\label{eqn:GFB1}
\end{equation}
where $\x^{(m)}_{N_m}$ denotes
$(x^{(m)}_1, \dots, x^{(m)}_{N_m})$, 
and (\ref{eqn:GFB1}) 
defines the {\it spatio-temporal correlation functions}
$\rho_{\xi}(\cdot)$ for the process $(\Xi(t), t \in \N_0, \P_{\xi})$.
Given an integral kernel
$$
\bK(s,x;t,y); 
\quad 
(s,x), (t,y) \in \N_0 \times \Z,
$$
the {\it Fredholm determinant} is defined as
\begin{eqnarray}
&& \mathop{{\rm Det}}_
{\substack{
(s,t)\in \{t_1, \dots, t_M\}, \\
(x,y)\in \Z^2}
}
 \Big[\delta_{st} \delta_x(\{y\})
+ \bK(s,x;t,y) \chi_{t}(y) \Big]
\nonumber\\
&& 
=\sum_
{\substack
{N_m \geq 0, \\ 1 \leq m \leq M} }
\sum_
{\substack
{\x^{(m)}_{N_m} \in \Z^{N_m} \cap \W_{N_m}, 
\\ 1 \leq m \leq M} }
\prod_{m=1}^{M}
\prod_{j=1}^{N_{m}} 
\chi_{t_m} \Big(x_{j}^{(m)} \Big)
\det_{\substack
{1 \leq j \leq N_{m}, 1 \leq k \leq N_{n}, \\
1 \leq m, n \leq M}
}
\Bigg[
\bK(t_m, x_{j}^{(m)}; t_n, x_{k}^{(n)} )
\Bigg].
\nonumber\\
\label{eqn:Fredholm1}
\end{eqnarray}

By the reducibility of determinantal martingales 
(Lemma \ref{thm:reducibility})
and a combinatorial argument, we can prove the 
following identity.
\begin{lem}
\label{thm:Fredholm}
Let $\u \in \Z_{\rm e}^N \cap \W_N$ and
$\xi=\sum_{j=1}^N \delta_{u_j}$.
For $M \in \N, t_m \in \N, 1 \leq m \leq M$,
$t_1 < \cdots < t_M \leq T \in \N$,
\begin{eqnarray}
&& \rE_{\u} \left[
\prod_{m=1}^M \prod_{j=1}^N
\{1+\chi_{t_m}(S_j(t_m)) \}
\cD_{\xi}(T, \S(T)) \right]
\nonumber\\
&& \qquad =
\mathop{{\rm Det}}_
{\substack{
(s,t)\in \{t_1, \dots, t_M\}, \\
(x,y)\in \Z^2}
}
 \Big[\delta_{st} \delta_x(y)
+ \mbK_{\xi}(s,x;t,y) \chi_{t}(y) \Big],
\nonumber
\end{eqnarray}
where $\mbK_{\xi}$ is given by (\ref{eqn:K1})
with (\ref{eqn:mart1b}).
\end{lem}
The same identity was proved for continuous-time DMR
in Section 2 of \cite{Kat14}.
So we omit the proof of Lemma \ref{thm:Fredholm}
for discrete-time DMR.

Now we arrive at one of the main theorems of the present paper.
\begin{thm} 
\label{thm:main1}
For any initial configuration $\xi \in \mM_0$
with $\xi(\Z_{\rm e})=N \in \N$, 
the noncolliding RW, $(\Xi(t), t \in \N_0, \P_{\xi})$
is determinantal with the correlation kernel
(\ref{eqn:K1}) with (\ref{eqn:mart1b}) 
in the sense that
the moment generating function (\ref{eqn:GF1})
is given by a Fredholm determinant
\begin{equation}
\Psi_{\xi}^{\t}[\f]
=\mathop{{\rm Det}}_
{\substack{
(s,t)\in \{t_1,t_2,\dots, t_M\}^2, \\
(x,y)\in \Z^2}
}
 \Big[\delta_{st} \delta_x(\{y\})
+ \mbK_{\xi}(s,x;t,y) \chi_{t}(y) \Big],
\label{eqn:Fred}
\end{equation}
and then all spatio-temporal correlation functions
are given by determinants as
\begin{eqnarray}
&& \rho_{\xi} \Big(t_1,\x^{(1)}_{N_1}; \dots;t_M,\x^{(M)}_{N_M} \Big) 
=
\left\{ \begin{array}{l}
\displaystyle{
\det_{\substack
{1 \leq j \leq N_{m}, 1 \leq k \leq N_{n}, \\
1 \leq m, n \leq M}
}
\Bigg[
\mbK_{\xi}(t_m, x_{j}^{(m)}; t_n, x_{k}^{(n)} )
\Bigg],
}
\cr
\qquad \mbox{if} \quad
\x^{(m)}_{N_m} \in \Z_{\rm e}^{N_m} \cap \W_{N_m}, 
t_m=\mbox{even}, \cr
\qquad  \mbox{or} \quad
\x^{(m)}_{N_m} \in \Z_{\rm o}^{N_m} \cap \W_{N_m}, 
t_m=\mbox{odd},  \quad 1 \leq m \leq M, \cr
0, \quad \mbox{otherwise},
\end{array} \right.
\nonumber\\
\label{eqn:main1}
\end{eqnarray}
$t_m \in \N, 1 \leq m \leq M$, 
$t_1 < \cdots < t_M$, and
$0 \leq N_m \leq N, 1 \leq m \leq M$.
\end{thm}
\noindent{\it Proof.} \,
By (\ref{eqn:Xi1}) with (\ref{eqn:Xi2}),
the moment generating function (\ref{eqn:GF1}) is 
written using (\ref{eqn:chi1}) as
$$
\Psi_{\xi}^{\t}[\f]
=\E_{\xi} \left[ \prod_{m=1}^N \prod_{j=1}^N
\{1+\chi_{t_m}(S^0_j(t_m)) \} \right].
$$
Proposition \ref{thm:DMR1} gives its DMR,
$$
\Psi_{\xi}^{\t}[\f]
=\rE_{\u} \left[
\prod_{m=1}^N \prod_{j=1}^N
\{1+\chi_{t_m}(S_j(t_m)) \} 
\cD_{\xi}(T, \S(T)) \right].
$$
Then Lemma \ref{thm:Fredholm} gives
(\ref{eqn:Fred}).
By definitions of correlation functions (\ref{eqn:GFB1})
and Fredholm determinant (\ref{eqn:Fredholm1}),
(\ref{eqn:main1}) is concluded from
(\ref{eqn:Fred}). The proof is completed. \qed

\SSC{Dynamics with an Infinite Number of Particles \label{sec:infinite}}
\subsection{Determinantal process with an infinite number of particles}

In this subsection, we will show that the noncolliding RW with an infinite number 
of particles can be well-defined as a determinantal process
for the initial configurations $\delta_{2 a \Z}, a \in \{2, 3, \dots\}$
given by (\ref{eqn:xi2a}). 
In order to that, we prepare infinite sequences of entire functions
and discrete-time martingales labeled by $k \in \Z$ below.

For a configuration $\xi =\sum_{j} \delta_{u_j} \in \mM_0$ we write
its restriction in $[-L, L] \subset \Z, L \in \N$ as
$\xi \cap [-L, L] \equiv \sum_{j: u_j \in [-L, L]} \delta_{u_j}$.
For each infinite-particle configuration (\ref{eqn:xi2a})
with $a \in \{2,3,\dots\}$, and $k \in \Z$,
a limit of the polynomial (\ref{eqn:Phi1})
\begin{equation}
\Phi_{\delta_{2a \Z}}^{2ak}(z)
\equiv \lim_{L \to \infty} \Phi_{\delta_{2a \Z} \cap [-L, L]}^{2a k}(z), 
\quad z \in \C
\label{eqn:hPhi2}
\end{equation}
exists and explicitly calculated as
\begin{eqnarray}
\Phi_{\delta_{2a \Z}}^{2ak}(z)
&=& \prod_{j \in \Z, j \not= k}
\frac{z-2aj}{2ak-2aj}
= \prod_{n \in \Z, n \not= 0}
\left(1 + \frac{z/2a -k}{n} \right)
\nonumber\\
&=& \frac{\sin (\pi(z/2a -k))}{\pi(z/2a-k)}
= \frac{1}{2 \pi} \int_{-\pi}^{\pi} d \lambda \,
e^{i \lambda(z/2a -k)},
\quad k \in \Z
\label{eqn:hPhi1}
\end{eqnarray}
by using the product formula of the sine function \cite{Lev96,KT10},
As the analytic continuation of (\ref{eqn:t_W2})
with respect to $\alpha$, 
\begin{equation}
\widetilde{\rE} [e^{-\lambda \widetilde{S}(t)}]
= \frac{1}{(\cos \lambda)^t},
\quad \lambda \in \left( -\frac{\pi}{2}, \frac{\pi}{2} \right),
\quad t \in \N_0, 
\label{eqn:t_W2b}
\end{equation}
implies that
$e^{i \lambda\{(y+i \widetilde{S}(t))/2a-k\}}$
is $d \lambda \times d \widetilde{\rP}$-integrable
for $a \geq 2$.
Then
\begin{equation}
\cM_{\delta_{2a \Z}}^{2ak}(t,y)
\equiv \widetilde{\rE} [
\Phi_{\delta_{2a \Z}}^{2ak} (y+i \widetilde{S}(t)) ],
\quad k \in \Z
\label{eqn:hM0}
\end{equation}
are well-defined and given by
\begin{equation}
\cM_{\delta_{2a \Z}}^{2ak}(t,y)
= \frac{1}{2 \pi} \int_{-\pi}^{\pi} d \lambda \,
\frac{e^{i \lambda(y/2a-k)}}{[\cos(\lambda/2a)]^t},
\quad k \in \Z.
\label{eqn:hM1}
\end{equation}
Since $|\cM_{\delta_{2a \Z}}^{2ak}(t,y)| \leq 2^{t/2}, a \in \{2,3, \dots\}$,
$|\sum_{j \in \Z} p(s,x|2aj) \cM_{\delta_{2a \Z}}^{2aj}(t,y)| < \infty$
for any $(s,t) \in \N^2$, $s, t \leq T \in \N$, $(x,y) \in \Z^2$.
Then the kernel
\begin{equation}
\mbK_{\delta_{2a\Z}}(s,x;t,y)
= \left\{ \begin{array}{l}
\displaystyle{
\sum_{j \in \Z} p(s, x|2a j) \cM_{\delta_{2a \Z}}^{2aj}(t,y)
- \1(s>t) p(s-t,x|y),
}
\cr
\qquad \qquad \qquad \mbox{if} \quad
(s,x), (t,y) \in \N_0 \times \Z, \quad s+x, t+y \in \Z_{\rm e},
\cr
0,
\qquad \qquad \quad \mbox{otherwise}, 
\end{array} \right.
\label{eqn:hK1}
\end{equation}
defines the moment generating function of the process by
the Fredholm determinant
$$
\Psi_{\delta_{2a\Z}}^{\t}[\f]
=\mathop{{\rm Det}}_
{\substack{
(s,t)\in \{t_1,t_2,\dots, t_M\}^2, \\
(x,y)\in \Z^2}
}
 \Big[\delta_{st} \delta_x(y)
+ \mbK_{\delta_{2a\Z}}(s,x;t,y) \chi_{t}(y) \Big]
$$
for any integer $M \in \N$,
a sequence of times
$\t=(t_1,\dots,t_M) \in \N^M$ with 
$t_1 < \cdots < t_M \leq T \in \N$,
and a sequence of bounded functions
$\f=(f_{t_1},\dots,f_{t_M})$
with (\ref{eqn:chi1}).
It implies that $\P_{\delta_{2a\Z}}$ is determined
in the sense of finite dimensional distributions.

\begin{prop}
\label{thm:inf_det}
For each $a \in \{2,3, \dots\}$, 
the noncolliding RW started at $\delta_{2a \Z}$,
denoted by $(\Xi(t), t \in \N_0, \P_{\delta_{2a\Z}})$,
is well-defined as a determinantal process
with the correlation kernel (\ref{eqn:hK1}).
\end{prop}

It is readily shown by Lemma \ref{thm:Fujita1}
(discrete It\^o's formula) that
if $(S(t))_{t \in \N_0}$ is a RW,
$(\cM_{\delta_{2a \Z}}^{2ak}(t, S(t)))_{t \in \N_0}, k \in \Z$
are discrete-time martingales, if $a \in \{2,3, \dots\}$.
Let $(S_j(t))_{t \in \N_0}$, $j \in \Z$ be an infinite
sequence of independent RWs.
Then we have an infinite sequence of independent
martingales with discrete time,
\begin{equation}
(\cM_{\delta_{2a \Z}}^{2ak}(t, S_j(t)))_{t \in \N_0}, 
\quad k \in \Z, 
\label{eqn:hmart1}
\end{equation}
for each $a \in \{2,3, \dots\}$ and $j \in \Z$.
We write the labeled configuration
$(2a j)_{j \in \Z}$ 
with an infinite number of particles as $2a \Z$,
and under $\rP_{2a \Z}$,
$S_j(0)=2aj, j \in \Z$. Then,
for any $t \in \N_0$, 
\begin{eqnarray}
\rE_{2a \Z} \Big[
\cM_{\delta_{2a \Z}}^{2ak} (t, S_j(t)) \Big]
&=& \rE_{2a \Z} \Big[
\cM_{\delta_{2a \Z}}^{2ak} (0, S_j(0)) \Big]
\nonumber\\
&=& 
\cM_{\delta_{2a \Z}}^{2ak} (0, 2a j) 
\nonumber\\
&=& \delta_{jk}, \quad j, k \in \Z.
\label{eqn:hM2}
\end{eqnarray}
Fix $N \in \N$. 
For $\J \subset \I_N$, define
the determinantal martingale of (\ref{eqn:hmart1})
\begin{equation}
\cD_{\delta_{2a\Z}}(t,\S_{\J}(t))
=\det_{j, k \in \J}
\Big[\cM_{\delta_{2a \Z}}^{k}(t, S_j(t)) \Big], 
\quad t \in \N_0.
\label{eqn:hDM1}
\end{equation} 
Let $t \in \N, t \leq T \in \N$,
$N' \in \N, N' < N$, 
and $F_{N'}$ be a symmetric bounded function on $Z^{N'}$.
Then the reducibility
\begin{eqnarray}
&& \sum_{\J \subset \I_N, |\J|=N'}
\rE_{2a \Z} \left[
F_{N'}(\S_{\J}(t))
\cD_{\delta_{2a \Z}}(T, \S_{N}(T)) \right]
\nonumber\\
&& \quad =
\sum_{\J \subset \I_N, |\J|=N'}
\rE_{2a \Z} \left[
F_{N'}(\S_{\J}(t))
\cD_{\delta_{2a \Z}}(T, \S_{\J}(T)) \right]
\nonumber\\
&& \quad
= \int_{\W_{N'}} \delta^{\otimes N'}_{2a\Z} (d\v)
\rE_{\v} \left[
F_{N'}(\S_{N'}(t))
\cD_{\delta_{2a\Z}}(T, \S_{N'}(T)) \right].
\label{eqn:reducibility2}
\end{eqnarray}
holds as well as Lemma \ref{thm:reducibility}.
Note that the last expression of (\ref{eqn:reducibility2})
does not change even if we replace $N$ in the 
first line of (\ref{eqn:reducibility2}) 
by any other integer $\widehat{N}$ with 
$\widehat{N} > N$.
Based on such consistency in reduction of DMRs
and the fact (\ref{eqn:hPhi2}), 
the DMR is valid also for 
the noncolliding RW with an infinite number of particles.

\vskip 0.3cm
\begin{prop}
\label{thm:inf_noncRW}
Assume that $F$ is represented as
$$
F(\Xi(\cdot))
=G \left( \sum_{x \in \Z} \phi_1(x) \Xi(t_1, x),
\dots, \sum_{x \in \Z} \phi_M(x) \Xi(t_M,x) \right),
$$
where $G$ is a polynomial function on $\R^M, M \in \N$
and $\phi_{m}, 1 \leq m \leq M$
are real-valued bounded functions with finite supports on $\Z$.
Then the expressions (\ref{eqn:DMR1}) are valid 
also in the cases with $\xi=\delta_{2 a \Z}$ and
$\u=2a \Z$,$a \in \{2,3,\dots\}$, 
even though $N= \delta_{2a\Z}(\Z)=\infty$.
\end{prop}
\vskip 0.3cm

\noindent
Proof is given in the similar way to that
given for Corollary 1.3 in \cite{KT13}.

\subsection{Relaxation to equilibrium dynamics}

Now we prove that the infinite-particle systems
$(\Xi(t), t \in \N_0, \P_{\delta_{2 a \Z}}), a \in \{2, 3, \dots\}$,
constructed in the previous subsection show 
relaxation phenomena to the equilibrium determinantal processes
with discrete analogues of the extended sine-kernel (\ref{eqn:K3b}).

Since the transition probability of RW (\ref{eqn:tp1})
is a unique solution of the difference equation
$$
p(t+1,y|x)=\frac{1}{2} 
[ p(t,y-1|x)+p(t,y+1|x)],
\quad t \in \N_0, \quad x, y \in \Z,
$$
with the initial condition $p(0, y|x)=\delta_{x y}$,
it has the following expressions,
\begin{eqnarray}
p(t,y|x) &=& \frac{1}{2 \pi}
\int_{-\pi}^{\pi} dk \,
e^{ik(y-x)} (\cos k)^t
\nonumber\\
&=& \frac{1}{4 \pi a}
\int_{-2a \pi}^{2 a \pi} d \theta \,
e^{i \theta(y-x)/2a}
\left[ \cos \left( \frac{\theta}{2a} \right) \right]^t
\nonumber\\
&=& \int_0^1 du \, \cos(u \pi(y-x)) [\cos(u \pi)]^{t},
\label{eqn:pB1}
\end{eqnarray}
where $a \in \N$.
Note that the integral representations (\ref{eqn:pB1})
of (\ref{eqn:tp1}) are valid for any 
$t \in \N_0, x, y \in \Z$.
Then combining with (\ref{eqn:hM1}) we have
$$
\sum_{j \in \Z} p(s,x|2aj) 
\cM_{\delta_{2a \Z}}^{2aj}(t,y)
= \frac{1}{8 \pi^2 a} \sum_{j \in \Z}
\int_{-2a \pi}^{2a \pi} d \theta 
\int_{-\pi}^{\pi} d \lambda \,
e^{i \theta(x/2a-j) + i \lambda(y/2a-j)}
\frac{[\cos (\theta/2a)]^s}{[\cos (\lambda/2a)]^t}.
$$
We rewrite the first line of (\ref{eqn:hK1})
as follows: for $(s,x), (t,y) \in \N_0 \times \Z$,
$s+x, t+y \in \Z_{\rm e}$,
\begin{equation}
\mbK_{\delta_{2a \Z}}(s,x;t,y)
+ \1(s>t)p(s-t, x|y)
= G(s,x;t,y)+R(s,x;t,y)
\label{eqn:K2a}
\end{equation}
with
$$
G(s,x;t,y)
= \frac{1}{4 \pi^2 a}
\int_{|\theta| \leq \pi} d \theta
\int_{|\lambda| \leq \pi} d \lambda \, 
\frac{e^{i(\theta x+ \lambda y)/2a}}
{[\cos(\lambda/2a)]^{t-s}}
\sum_{j \in \Z} e^{-i(\theta+\lambda)j}
\left[ \frac{\cos(\theta/2a)}{\cos(\lambda/2a)}\right]^s,
$$
and
$$
R(s,x;t,y)
= \frac{1}{8 \pi^2 a} \sum_{j \in \Z}
\int_{\pi <|\theta| < (2a-1) \pi} d \theta
\int_{|\lambda| \leq \pi} d \lambda \, 
\frac{e^{i(\theta x+ \lambda y)/2a}}
{[\cos(\lambda/2a)]^{t-s}}
e^{-i(\theta+\lambda)j}
\left[ \frac{\cos(\theta/2a)}{\cos(\lambda/2a)}\right]^s.
$$
Since 
$\sum_{j \in \Z} e^{-i(\theta+\lambda)j}
=2\pi \delta_{-\lambda}(\{\theta \})$
for $\theta, \lambda \in (-\pi, \pi]$, we obtain
\begin{equation}
G(s,x; t,y)
= \frac{1}{2 \pi a}
\int_{-\pi}^{\pi} d \lambda \,
\frac{e^{i \lambda(y-x)/2a}}
{[\cos (\lambda/2a)]^{t-s}}
\equiv \cG(t-s,y-x).
\label{eqn:K2d}
\end{equation}
On the other hand, when 
$\pi < |\theta| < (2a-1) \pi$
and $|\lambda| \leq \pi$,
$
|\cos(\theta/2a)/\cos(\lambda/2a)| < 1.
$
Then for any fixed $s,t \in \N$,
$
|R(s+n, x; t+ n, y)| \to 0$
as $n \to \infty$ 
uniformly on any $(x,y) \in \Z^2$
and it implies
\begin{equation}
\mbK_{\delta_{2a \Z}}(s+n, x; t+n, y)
\to \bK_{\rho}(t-s, y-x)
\quad \mbox{as $n \to \infty$},
\label{eqn:limit1}
\end{equation}
where
\begin{eqnarray}
\bK_{\rho}(t-s, y-x)
&=&\cG(t-s,y-x)-\1(s>t) p(s-t, x|y)
\nonumber\\
&=& 
2 \int_{0}^{\rho} d u \,
\frac{\cos(\pi u (y-x))}{[\cos (\pi u)]^{t-s}}
-\1(s>t) p(s-t,x|y), 
\label{eqn:K3a}
\end{eqnarray}
if $s+x, t+y \in \Z_{\rm e}$, and
$\bK_{\rho}(t-s,y-x)=0$, otherwise, 
with the density on $\Z$, 
\begin{equation}
\rho=\frac{1}{2 a}, \quad a \in \{2,3, \dots\}.
\label{eqn:rho1}
\end{equation}
By (\ref{eqn:pB1}) and others, 
we can see that (\ref{eqn:K3a}) is written as (\ref{eqn:K3b}).

The convergence of the correlation kernel 
(\ref{eqn:limit1}) implies the convergence of
generating function for correlation functions
$\Psi_{\delta_{2a \Z}}^{\t}[\f]$,
and thus the convergence of the determinantal
process to an equilibrium determinantal process.
This is an example of {\it relaxation phenomena}
\cite{KT09,KT10,KT11,Esa14,Kat14b}.

In order to state the result, we define 
determinantal point processes on $\Z$.

\begin{df}
\label{thm:sine}
Let $\sharp={\rm e}$ or ${\rm o}$.
For a given density $0 < \rho < 1/2$, 
the probability measures
$\mu^{\rm sin}_{\rho, \sharp}$ on $\Z$
are defined as determinantal point processes
with the sine kernels
\begin{equation}
\bK^{\rm sin}_{\rho, \sharp}(y-x)
=\left\{
\begin{array}{ll}
\displaystyle{\frac{2 \sin(\pi \rho(y-x))}{\pi (y-x)}},
\quad & \mbox{if $x,y \in \Z_{\sharp}$},
\cr
0, \quad & \mbox{otherwize}.
\end{array}
\right.
\label{eqn:sine1}
\end{equation}
\end{df}

\begin{thm}
\label{thm:equilibrium}
For each $a \in \{2,3, \dots\}$, 
the process $(\Xi(t), t \in \N_0, \P_{\delta_{2a \Z}})$
shows a relaxation phenomenon to equilibrium state
such that
\begin{eqnarray}
\Xi(2n) &\Rightarrow& \mu^{\rm sin}_{\rho, {\rm e}},
\nonumber\\
\Xi(2n+1) &\Rightarrow& \mu^{\rm sin}_{\rho, {\rm o}},
\quad \mbox{as $n \to \infty$}
\nonumber
\end{eqnarray}
with $\rho=1/2a$.
The equilibrium process, denoted by $(\Xi(t), t \in \Z, \bP_{\rho})$,  
is time-reversible and 
determinantal with the correlation kernel given by (\ref{eqn:K3b}).
\end{thm}
\vskip 0.3cm

Here we note that the local densities of particles
(the one-point correlation functions) 
in $\mu^{\sin}_{\rho, \sharp}$ and
in $(\Xi(t), t \in \Z, \bP_{\rho})$ are obtained
from the expressions (\ref{eqn:sine1}) and (\ref{eqn:K3b}) for
correlation kernels, respectively, by taking the limits as
\begin{eqnarray}
\mu^{\rm sin}_{\rho, \sharp}(\{x\})
&=& \lim_{y \to x} \bK^{\rm sin}_{\rho, \sharp}(y-x)
=\left\{ \begin{array}{ll}
2 \rho, \qquad \mbox{if $x \in \Z_{\sharp}$}, \cr
0, \qquad \mbox{otherwize},
\end{array} \right.
\quad \mbox{$\sharp={\rm e}$ or ${\rm o}$},
\nonumber\\
\bP_{\rho}[\Xi(s,\{x\})=1]
&=& \lim_{\substack{t \to s, \cr y \to x}} \bK_{\rho}(y-x, t-s)
=\left\{ \begin{array}{ll}
2 \rho, \qquad \mbox{if $s+x \in \Z_{\rm e}$}, \cr
0, \qquad \mbox{otherwize}.
\end{array} \right.
\nonumber
\end{eqnarray}
On the spatio-temporal plane $(t,x) \in \Z^2$,
the equilibrium state makes a homogeneous
bipartite lattice.

\SSC{Convergence to the Dyson Model}\label{sec:convergence}

In this final section, we will discuss the convergence of
noncolliding RWs to the continuous version ({\it i.e.} the Dyson model)
in the sense of Donsker's invariant principle from the viewpoint of DMR.

For $n \in \N$, define scaled discrete-processes as
\begin{eqnarray}
&& S_j^{(n)}(t) = \frac{1}{n}S_j(n^2 t), \quad
\widetilde{S}_j^{(n)}(t) = \frac{1}{n}\widetilde{S}_j(n^2 t), 
\nonumber\\
&& Z_j^{(n)}(t) = S_j^{(n)}(t)+ i \widetilde{S}_j^{(n)}(t),
\quad t \in \N_0, \quad 1 \leq j \leq N.
\label{eqn:appr1}
\end{eqnarray}
We set $S_j(0)=n u_j,1 \leq j \leq N$.
Since
\begin{eqnarray}
&& \rE[\zeta_j(t)]= \widetilde{\rE}[\widetilde{\zeta}_j(t)]=0,
\nonumber\\
&& \rE[\zeta_j(t)^2]= \widetilde{\rE}[\widetilde{\zeta}_j(t)^2]=1,
\quad t \in \N, \quad 1 \leq j \leq N,
\nonumber
\end{eqnarray}
{\it Donsker's invariance principle} \cite{Bil68,RY05} proves
the convergence in distribution 
\begin{equation}
S_j^{(n)}(\cdot) \dto B_j(\cdot), \quad
\widetilde{S}_j^{(n)}(\cdot) \dto \widetilde{B}_j(\cdot), \quad
Z_j^{(n)}(\cdot) \dto \cB_j(\cdot),
\quad \mbox{as $n \to \infty$},
\label{eqn:Donsker1}
\end{equation}
where $B_j(\cdot)$ and $\widetilde{B}_j(\cdot)$
are independent BMs with
$B_j(0)=u_j, \widetilde{B}_j(0)=0, 1 \leq j \leq N$,
and $\cB_j$ denotes the complex BMs,
$\cB_j(\cdot)=B_j(\cdot)+i \widetilde{B}_j(\cdot)$, $1 \leq j \leq N$.
For $\Phi_{\xi}^{u_k}(\cdot), 1 \leq k \leq N$
are polynomials and thus continuous functions,
(\ref{eqn:Donsker1}) implies
\begin{equation}
\Phi_{\xi}^{u_k}
(Z_j^{(n)}(\cdot)) \dto
\Phi_{\xi}^{u_k}(\cB_j(\cdot))
\quad \mbox{as $n \to \infty$},
\quad 1 \leq j, k \leq N.
\label{eqn:Donsker2}
\end{equation}

For each $n \in \N$, let 
$\S^{0 \, (n)}(\cdot)=(S_1^{0 \, (n)}(\cdot), \dots, S_N^{0 \,(n)}(\cdot))$
be the $N$-particle scaled RW 
{\it conditioned never to collide with each other}
started at $\u=(u_1, \dots, u_N) \in \Z_{\rm e}^N \cap \W_N$
and put 
$\Xi^{(n)}(t, \cdot)=\sum_{j=1}^N \delta_{S_j^{0 \, (n)}(t)}(\cdot)$,
$t \in \N_0$.
Then we have a series of scaled noncolliding RWs,
$(\Xi^{(n)}(t), t \in\N_0, \P_{\xi})$, $n \in \N$,
each of which has DMR
\begin{eqnarray}
\E_{\xi} \left[ F \left(\Xi^{(n)}(\cdot) \right) \right]
&=& \rE_{n \u} \left[F \left( \sum_{j=1}^{N} \delta_{S_j^{(n)}(\cdot)} \right)
\cD_{\xi}( n^2 T, \S^{(n)}(T)) \right]
\nonumber\\
&=& \bE_{n \u} \left[F \left( \sum_{j=1}^{N} \delta_{\Re Z_j^{(n)}(\cdot)} \right)
\det_{1 \leq j, k \leq N}
[\Phi_{\xi}^{u_k}(Z_j^{(n)}(T))] \right]
\label{eqn:DMRn1}
\end{eqnarray}
for any $\cF(t)$-measurable bounded function $F$
for any $t \in \N, t \leq T \in \N$.
Let $(\vXi(t), t \in [0, \infty), \cP_{\xi})$
be the Dyson model
started at $\xi=\sum_{j=1}^N \delta_{u_j} \in \mM_0$
with $\u=(u_1, \dots, u_N) \in \Z_{\rm e}^N \cap \W_N$.
That is, 
$\vXi(t, \cdot)=\sum_{j=1}^N \delta_{X_j(t)}(\cdot), t \in [0, \infty)$,
where $\X(\cdot)=(X_1(\cdot), \dots, X_N(\cdot))$
is a unique solution of the SDEs (\ref{eqn:Dyson1}) 
under the initial configuration $\X(0)=\u \in \Z_{\rm e}^N \cap \W_N$.
By the invariance principle (\ref{eqn:Donsker1}), (\ref{eqn:Donsker2}),
if $F$ is continuous, 
the DMRs given  by the RHS of (\ref{eqn:DMRn1})
converge to the {\it complex BM representation}
for $(\vXi(t), t \in [0, \infty), \cP_{\xi})$
given by Theorem 1.1 in \cite{KT13}.
Since the complex BM representation is a special
case of DMR (see Remark 4 and a comment mentioned
just after Proposition \ref{thm:DMR1}), 
we will say that
\begin{equation}
\mbox{$(\Xi^{(n)}(t), t \in \N_0, \P_{\xi})$
converges to $(\vXi(t), t \in [0, \infty), \cP_{\xi})$ in DMR.}
\label{eqn:conv1}
\end{equation}

As shown in Section \ref{sec:det_prop},
the DMR gives a Fredholm determinantal expression
for any generating function of multitime correlation functions.
Then (\ref{eqn:conv1}) implies the convergence
in the sense of finite dimensional distributions.
It also implies the convergence as determinantal processes.
By the convergence of processes (\ref{eqn:Donsker1}), 
the following convergence of functions are concluded;
if $p(n^2 t, n y|nx) \not=0$,
$$
\rP_{n x}[S^{(n)}(t) \in dy]
= p(n^2 t, ny|n x) n dy \to p^{\rm BM}(t,y|x),
$$
and
\begin{eqnarray}
&& \widetilde{\rm E} \left[
\Phi_{\xi}^{u_k}(x+i \widetilde{S}_j^{(n)}(t)) \right]
= \cM_{\xi}^{u_k}(n^2 t, x)
\nonumber\\
&& \qquad \qquad \to
\int_{\R} dv \, p^{\rm BM}(t, v|0)
\Phi_{\xi}^{u_k}(x+iv)
\equiv \sfM_{\xi}^{u_k}(t, x)
\nonumber
\end{eqnarray}
as $n \to \infty$ with (\ref{eqn:pBM}).
Therefore, the correlation kernel of 
the Dyson model, $(\vXi(t), t \in [0, \infty), \cP_{\xi})$,
is determined as the limit
of the kernels of $(\Xi^{(n)}(t), t \in \N_0, \P_{\xi})$
of the form (\ref{eqn:K1}), 
\begin{equation}
\sfK_{\xi}(s,x;t,y)
=\sum_{j=1}^N
p^{\rm BM}(s,x|u_j) \sfM_{\xi}^{u_j}(t,y)
- \1(s>t) p^{\rm BM}(s-t, x|y),
\label{eqn:KDyson2}
\end{equation}
$(s,x), (t,y) \in [0, \infty) \times \R$.
The limit (\ref{eqn:KDyson2}) is exactly the same
as the correlation kernel of the Dyson model 
given as Eq.(2.2) in \cite{KT10}
for general $\xi \in \mM_0$,
$\xi(\R)=N \in \N$,
which was obtained by using the multiple Hermite polynomials.

As claimed by Proposition \ref{thm:inf_noncRW},
DMR is valid for $(\Xi(t), t \in \N_0, \P_{\delta_{2 a \Z}})$,
$a \in \{2,3,\dots\}$.
Then we will conclude that
\begin{equation}
\mbox{$(\Xi^{(n)}(t), t \in \N_0, \P_{\delta_{2 a \Z}})$
converges to $(\vXi(t), t \in [0, \infty), \cP_{\delta_{2 a \Z}})$ in DMR},
\label{eqn:convB1}
\end{equation}
if $a \in \{2,3, \dots\}$.
From (\ref{eqn:hM0}) with (\ref{eqn:hPhi1}),
\begin{eqnarray}
&&\widetilde{\rE} \left[
\Phi_{\delta_{2 a \Z}}^{2 a k}
(x+i \widetilde{S}^{(n)}(t)) \right]
= \cM_{\delta_{2 a \Z}}^{2 a k}
(n^2 t, x) 
\nonumber\\
&& \qquad \to
\int_{\R} dv \, p^{\rm BM}(t, v|0)
\Phi_{\delta_{2 a \Z}}^{2 a k}(x+iv)
\equiv \sfM_{\delta_{2 a \Z}}^{2 a k}(t, x)
\label{eqn:convB2}
\end{eqnarray}
as $n \to \infty$.
We find
$$
\sfM_{\delta_{2 a \Z}}^{2ak}(t,y)
= \frac{1}{2 \pi} \int_{-\pi}^{\pi} d \lambda \,
e^{\lambda^2 t/8 a^2 + i \lambda(y/2a-k)},
\quad (t,y) \in [0, \infty) \times \R,
\quad k \in \Z.
$$
Then the Dyson model with an infinite number of particles,
$(\vXi(t), t \in [0, \infty), \cP_{\delta_{2 a \Z}})$,
is determinantal and its correlation kernel
is determined as
\begin{equation}
\sfK_{\delta_{2 a \Z}}(s,x;t,y)
= \sum_{j \in \Z} p^{\rm BM}(s,x|2aj)
\sfM_{\delta_{2 a \Z}}^{2 a k}(t, y)
- \1(s>t) p^{\rm BM}(s-t, x|y),
\label{eqn:Kinf1}
\end{equation}
$(s,x), (t,y) \in [0, \infty) \times \R$.
Let
$$
\vartheta_3(v, \tau)
= \sum_{j \in \Z} e^{2 \pi i v j+ \pi i \tau j^2},
\quad \Im \tau > 0, 
$$
which is a version of the Jacobi theta function.
If we use the reciprocity relation
$$
\vartheta_3(v,\tau)= 
\vartheta_3 \left( \frac{v}{\tau}, - \frac{1}{\tau} \right)
e^{-\pi i v^2/\tau} \sqrt{ \frac{i}{\tau}}
$$
(see, for example, Section 10.12 in \cite{AAR99}),
we can obtain the expression
\begin{equation}
\sum_{j \in \Z} p^{\rm BM}(s, x|2aj)
\sfM_{\delta_{2a \Z}}^{2a j} (t, y)
= \frac{\rho}{2 \pi}
\int_{-\pi}^{\pi} d \lambda \,
e^{\lambda^2 \rho^2 (t-s)/2 + i \lambda \rho (y-x)}
\vartheta_3( \rho x- i \lambda \rho^2 s, 2 \pi i \rho^2 s),
\nonumber
\end{equation}
where $\rho$ is the density of particles given by (\ref{eqn:rho1}).
Then (\ref{eqn:Kinf1}) is written as
\begin{eqnarray}
&& \sfK_{\delta_{2 a \Z}}(s,x;t,y)
= \sfK_{\rho}^{\rm sin}(t-s, y-x)
\nonumber\\
&& \qquad \qquad 
+ \frac{\rho}{2 \pi} \int_{-\pi}^{\pi} d \lambda \,
e^{\lambda^2 \rho^2 (t-s)/2 + i \lambda \rho(y-x)}
\{ \vartheta_3(\rho x- i \lambda \rho^2 s, 2 \pi i \rho^2 s)
-1 \},
\label{eqn:Kinf2}
\end{eqnarray}
$(s,x), (t,x) \in [0, \infty) \times \R$,
where
\begin{equation}
\sfK_{\rho}^{\rm sin}(t-s, y-x)
= \left\{ \begin{array}{ll}
\displaystyle{ \int_0^{\rho} du \,
e^{\pi^2 u^2(t-s)/2} \cos(\pi u(y-x))},
& \mbox{if $s<t$}, \cr
& \cr
\displaystyle{ \frac{\sin(\pi \rho(y-x))}{\pi(y-x)}},
& \mbox{if $s=t$},\cr
& \cr
\displaystyle{ -\int_{\rho}^{\infty} du \,
e^{\pi^2 u^2(t-s)/2} \cos(\pi u \pi(y-x))},
& \mbox{if $s>t$}. \cr
\end{array} \right.
\label{eqn:Kinf3}
\end{equation}
The correlation kernel (\ref{eqn:Kinf2}) coincide
with Eq.(1.5) in \cite{KT10} if we set $\rho=1$.
The kernel (\ref{eqn:Kinf3}) is called
the {\it extended sine kernel} with density $\rho$ 
(see Section 11.7.1 in \cite{For10}),
which is a continuum limit of (\ref{eqn:K3b}).
The relaxation phenomenon associated with
$
\lim_{\tau \to \infty}
\sfK_{\delta_{2 a Z}}(s+\tau, x; t + \tau, y)
= \sfK_{1}^{\rm sin}(t-s, y-x)
$
was studied in \cite{KT10}.

The above shows that the convergence in DMR implies
the convergence in the sense of finite dimensional distributions
and that as determinantal processes.
As demonstrated by Proposition 1.4 and Theorem 1.5 in \cite{KT13},
DMR is useful to test the Kolmogorov criterion for tightness.
Relations between the present convergence in DMR 
and the previous results concerning convergence to the Dyson model
\cite{KT03,BS07,PAT07} will be discussed elsewhere.

\vskip 0.5cm
\noindent{\bf Acknowledgements} \quad
The present author would like to thank 
T. Shirai, H. Osada, H. Tanemura, and S. Esaki
for useful discussions.
This work is supported in part by
the Grant-in-Aid for Scientific Research (C)
(No.21540397 and No.26400405) of Japan Society for
the Promotion of Science.



\begin{thebibliography}{99}
\bibitem{AAR99}
Andrews, G. E., Askey, R., Roy, R.:
Special functions.  
Cambridge University Press, Cambridge, U.K. (1999)

\bibitem{Bai00}
Baik, J.: 
Random vicious walks and random matrices. 
Commun. Pure Appl. Math. {\bf 53}, 1385-1410 (2000)

\bibitem{BS07}
Baik, J., Suidan, T. M.:
Random matrix central limit theorems for
nonintersecting random walks.
Ann. Probab. {\bf 35}, 1807-1834 (2007)

\bibitem{BKPV09}
Ben Hough, J., Krishnapur, M., Peres, Y., Vir\'ag, B.:
Zeros of Gaussian Analytic Functions and
Determinantal Point Processes.
Amer. Math. Soc., Providence (2009)

\bibitem{Bil68}
Billingsley, P.:
Convergence of Probability Measures. 
Wiley, New York (1968)

\bibitem{CK03}
Cardy, J., Katori, M.:
Families of vicious walkers. 
J. Phys. A {\bf 36}, 609-629 (2003)

\bibitem{Dys62}
Dyson, F. J. :
A Brownian-motion model for the eigenvalues of a random matrix.
J. Math. Phys. {\bf 3}, 1191-1198 (1962)

\bibitem{EK08}
Eichelsbacher, P., K\"onig, W.:
Ordered random walks.
Electron. J. Probab. {\bf 13}, no.46, 1307-1336 (2008)

\bibitem{Erd55}
Erd\'elyi, A. (editor) :
Higher Transcendental Functions. Bateman Manuscript Project,
vol. 3, McGraw-Hill, New York (1955)

\bibitem{Esa14}
Esaki, S.: Noncolliding system of continuous-time random walks.
Pacific J. Math. Industry {\bf 6}, 11/1-10 (2014)

\bibitem{EM98}
Eynard, B., Mehta, M. L. :
Matrices coupled in a chain: I.
Eigenvalue correlations.
J. Phys. A {\bf 31}, 4449-4456 (1998)

\bibitem{Fei12}
Feierl, T.:
The height of watermelons with wall. 
J. Phys. A {\bf 45}, 095003 (2012)

\bibitem{Fel66}
Feller, W.:
An Introduction to Probability Theory and Its Applications.
2nd edn. Volume II, Wiley, New York (1966)

\bibitem{Fis84}
Fisher, M.E.:
Walks, walls, wetting, and melting. 
J. Stat. Phys. {\bf 34}, 667-729, (1984)

\bibitem{For10}
Forrester, P. J.:
Log-gases and Random Matrices.
London Mathematical Society Monographs, 
Princeton University Press, Princeton (2010)

\bibitem{Fuj02}
Fujita, T. :
Stochastic Calculus for Finance
(in Japanese).
Kodansha, Tokyo (2002)

\bibitem{Fuj08}
Fujita, T., Kawanishi, Y.:
A proof of Ito's formula using a discrete Ito's formula.
Stud. Sci. Math. Hungr. {\bf 45}, 125-134 (2008)

\bibitem{Gra99}
Grabiner, D. J.:
Brownian motion in a Weyl chamber.
non-colliding particles, and random matrices. 
Ann. Inst. Henri Poincar\'e,
Probab. Stat. {\bf 35}, 177-204 (1999)

\bibitem{GM13}
Graczyk, P., Ma{\l}ecki, J.:
Multidimensional Yamada-Watanabe theorem 
and its applications to particle systems. 
J. Math. Phys. {\bf 54}, 021503/1-15 (2013)

\bibitem{HH68}
Harkness, W. L., Harkness, M. L.:
Generalized hyperbolic secant distributions.
J. Amer. Statist. Assoc. {\bf 63}, 329-337 (1968)

\bibitem{Joh01}
Johansson, K.:
Universality of the local spacing distribution
in certain ensembles of Hermitian Wigner matrices.
Commun. Math. Phys. {\bf 215}, 683-705 (2001)

\bibitem{Joh02}
Johansson, K.:
Non-intersecting paths, random tilings and random matrices. 
Probab. Theory Relat. Fields {\bf 123}, 225-280 (2002)

\bibitem{Joh05}
Johansson, K.:
Non-intersecting, simple, symmetric random walks
and the extended Hahn kernel.
Ann. Inst. Fourier {\bf 55}, 2129-2145 (2005)

\bibitem{Kat14}
Katori, M.:
Determinantal martingales and noncolliding diffusion processes.
Stochastic Process. Appl. {\bf 124}, 3724-3768 (2014)

\bibitem{Kat14b}
Katori, M.:
Elliptic determinantal process of type A.
Probab. Theory Relat. Fields, 
DOI 10.1007/s00440-014-0581-9

\bibitem{KT03}
Katori, M., Tanemura, H.:
Functional central limit theorems for vicious walkers.
Stoch. Stoch. Rep. {\bf 75}, 369-390 (2003)

\bibitem{KT04}
Katori, M., Tanemura, H.:
Symmetry of matrix-valued stochastic processes and
noncolliding diffusion particle systems. 
J. Math. Phys. {\bf 45}, 3058-3085 (2004)

\bibitem{KT07}
Katori, M. Tanemura, H.:
Noncolliding Brownian motion and determinantal processes.
J. Stat. Phys. {\bf 129}, 1233-1277 (2007)

\bibitem{KT09}
Katori, M., Tanemura, H.:
Zeros of Airy function and relaxation process.
J. Stat. Phys. {\bf 136}, 1177-1204 (2009)

\bibitem{KT10}
Katori, M., Tanemura, H.:
Non-equilibrium dynamics of Dyson's model
with an infinite number of particles.
Commun. Math. Phys. {\bf 293}, 469-497 (2010)

\bibitem{KT11}
Katori, M., Tanemura, H.:
Noncolliding squared Bessel processes.
J. Stat. Phys. {\bf 142}, 592-615 (2011)

\bibitem{KT11b}
Katori, M. and Tanemura, H.:
Noncolliding processes, matrix-valued processes
and determinantal processes. 
Sugaku Expositions (AMS) 
{\bf 24}, 263-289 (2011)

\bibitem{KT13}
Katori, M., Tanemura, H.:
Complex Brownian motion representation
of the Dyson model.
Electron. Commun. Probab. {\bf 18}, no.4, 1-16 (2013)

\bibitem{Koe05}
K\"onig, W.:
Orthogonal polynomial ensembles in probability theory.
Probab. Surveys {\bf 2}, 385-447 (2005)

\bibitem{KOR02}
K\"onig, W., O'Connell, N., Roch, S.:
Non-colliding random walks, tandem queues, and
discrete orthogonal polynomial ensembles.
Electron. J. Probab. {\bf 7} (1), 1-24 (2002)

\bibitem{Kra06}
Krattenthaler, C.:
Watermelon configurations with wall interaction :
exact and asymptotic results.
J. Phys. Conf. Series {\bf 42}, 179-212 (2006)

\bibitem{Lev96}
Levin, B. Ya.:
Lectures on Entire Functions. 
Translations of Mathematical Monographs, {\bf 150},
Amer. Math. Soc., Providence (1996)

\bibitem{Meh04}
Mehta, M. L.:
Random Matrices. 3rd edn. 
Elsevier, Amsterdam (2004)

\bibitem{NF98}
Nagao, T., Forrester, P.:
Multilevel dynamical correlation functions for Dyson's
Brownian motion model of random matrices.
Phys. Lett. {\bf A247}, 42-46 (1998)

\bibitem{NF02}
Nagao, T., Forrester, P.J.:
Vicious random walkers and a discretization of
Gaussian random matrix ensembles.
Nucl. Phys. B {\bf 620} [FS], 551-565 (2002)

\bibitem{Osa12}
Osada, H.:
Infinite-dimensional stochastic differential equations
related to random matrices.
Probab. Theory Relat. Fields {\bf 153}, 471-509 (2012)

\bibitem{Osa13}
Osada, H.:
Interacting Brownian motions in infinite dimensions
with logarithmic potentials.
Ann. Probab. {\bf 41}, 1-49 (2013)

\bibitem{Osa13b}
Osada, H.:
Interacting Brownian motions in infinite dimensions
with logarithmic interaction potentials II:
Airy random point field.
Stochastic Process. Appl. {\bf 123}, 813-838 (2013)

\bibitem{PAT07}
P\'erez-Abreu, V., Tudor, C.:
Functional limit theorems for trace processes in a Dyson Brownian motion.
Commun. Stoch. Anal. {\bf 1} (3), 415-428 (2007)

\bibitem{Pri09}
Privault, N.:
Stochastic Analysis in Discrete and Continuous Settings.
Lecture Notes in Mathematics {\bf 1982},
Springer, Berlin (2009)

\bibitem{RY05}
Revuz, D., Yor, M.:
Continuous Martingales and Brownian Motion. 3rd edn. 
Springer, New York (2005)

\bibitem{Sat99}
Sato, K.:
L\'evy Processes and Infinitely Divisible Distributions.
Cambridge University Press, Cambridge, U.K. (1999)

\bibitem{Sch00}
Schoutens, W.:
Stochastic Processes and Orthogonal Polynomials. 
Lecture Notes in Statistics {\bf 146}, 
Springer, New York (2000)

\bibitem{ST03}
Shirai, T., Takahashi, Y.: 
Random point fields associated with certain
Fredholm determinants I:
fermion, Poisson and boson point process.
J. Funct. Anal.
{\bf 205}, 414-463 (2003)

\bibitem{Sos00}
Soshnikov, A. : 
Determinantal random point fields.
Russian Math. Surveys {\bf 55}, 923-975 (2000)

\bibitem{Spo87}
Spohn, H.:
Interacting Brownian particles:
a study of Dyson's model.
In: Papanicolaou, G. (ed)
Hydrodynamic Behavior and Interacting Particle Systems, 
IMA Volumes in Mathematics and its Applications, vol.9, 
pp.151-179, Springer, Berlin (1987)

\bibitem{Tao12}
Tao, T.:
Topics in Random Matrix Theory. 
Amer. Math. Soc., Providence (2012)

\end{thebibliography}
\end{document}